\documentclass[10pt,a4paper]{article}
\makeatletter
\usepackage[dvipsnames]{xcolor}
\usepackage{lipsum}
\usepackage{tipa}
\usepackage{tikz}
\usepackage{tikz-cd}
\usetikzlibrary{matrix,arrows,decorations.pathmorphing}
\newcounter{magicrownumbers}

\usepackage[left=3cm,top=3cm,bottom=3cm,right=3cm,head=1cm,foot=1cm]{geometry}
\usepackage[pdftex,bookmarks=true,bookmarksnumbered=true]{hyperref}
\hypersetup{linkbordercolor=ForestGreen!25, citebordercolor=Goldenrod!60}
\DeclareMathAlphabet{\mathpzc}{OT1}{pzc}{m}{it}
\usepackage{graphicx}
\usepackage{stmaryrd}
\usepackage{twoopt}
\usepackage{amssymb}
\usepackage{amsmath}
\usepackage{amsthm}
\usepackage{mathabx} 
\usepackage{amsfonts}
\usepackage{amscd}
\usepackage{url}
\usepackage[all]{xy}

\newcommand{\5}{\hspace{0,5cm}}
\newcommand{\3}{\vspace{0,3cm}}

\newcommand{\la}{\langle}
\newcommand{\ra}{\rangle}

\newcommand{\fp}[2]{\mspace{1mu} {}_{#1} \mspace{-5mu} \times_{#2}}
\newcommand{\tx}[1]{\mathpzc{#1}}

\newcommand{\isonelio}[8]{\begin{tikzpicture}[scale=1.25 ] 
\path (3,1) node (b) {${#2}$} (0,1) node (a) {${#1}$} (3,0) node (d) {${#4}$} (0,0) node (c) {${#3}$};
\begin{scope}
\draw[->] (a) -- node [above] {${#5}$} (b);
\draw[->] (a) -- node [left] {${#7}$} (c) ;
\draw[->] (b) -- node [right] {${#8}$} (d) ;
\draw[->] (c) -- node [below] {${#6}$} (d); 
\end{scope}
\end{tikzpicture}}

\newcommand{\gpdhom}[7]{\begin{tikzpicture}
\path (3,1) node (a) {${#1}$} (0,0) node (b) {${#2}$} (2,0) node (c) {${#3}$} (4,0) node (d) {${#4}$} (6,0) node (e) {${#5}$}; 
\begin{scope}
\draw[->] (a) -- node [left] {${#6}$} (c) ; 
\draw[->] (a) -- node [right] {${#7}$} (d); 
\draw[->] (0.4,0.05) -- (1.6,0.05) ;
\draw[->] (0.4,-0.05) -- (1.6,-0.05) ;
\draw[->] (5.6,0.05) -- (4.4,0.05) ;
\draw[->] (5.6,-0.05) -- (4.4,-0.05) ;
\draw[->] (0,0.4) arc (180:90:64pt and 18pt);
\draw[->] (6,0.4) arc (0:90:64pt and 18pt);
\end{scope}
\end{tikzpicture}}

\title{On Noncommutative Geometry of Orbifolds}
\date{}
\author{Antti J. Harju\footnote{Copernicus Center for Interdisciplinary Studies, Krakow, Poland} \footnote{harjuaj@gmail.com}}

\begin{document}

\maketitle

\begin{abstract}
An orbifold is a Morita equivalence class of a proper {\' e}tale Lie groupoid. A unitary equivalence class of spectral triples over the algebra of smooth invariant functions are associated with any compact spin orbifold. In the case of an effective spin orbifold we construct a collection of spectral triples over the smooth convolution algebras of the representatives of the Morita equivalence class.\3

\noindent MSC 58B34, 22A22, 57R18  
\end{abstract}

\section*{Introduction}

\noindent \textbf{I.1.} The theory of spectral triples has proven to be a successful model for extensions of the theory of riemannian manifolds to the noncommutative geometric realm. A typical noncommutative geometric space is constructed by means of an analytic deformation of a function algebra on a manifold. Recently several examples of analytic deformations of singular spaces, such as orbifolds, have appeared in the literature, see e.g. \cite{A1}, \cite{A2}, \cite{BF12}, \cite{A3}. Spectral geometry has been discussed in the references, \cite{Har14b}, \cite{SV13}. The study of singular noncommutative deformations is somewhat out of the scope of the spectral triple formulation of noncommutative geometry since the axioms of the spectral triples are specific to smooth manifolds without singularities, \cite{Con94}, \cite{Con13}. The goal of this project is to develop Dirac spectral triples over function algebras of orbifolds in the differential geometric context. 

An orbifold (groupoid) is a Morita equivalence class of a proper {\' e}tale Lie groupoid, \cite{Moe02}, \cite{MM03}. Roughly speaking, a Morita equivalence class consist of those Lie groupoids which share the same orbit space, including its singularities, and therefore this geometric context is suitable for our purposes. All the singularities of an orbifold are of the finite type: the isotropy groups on the orbifold base are all finite. The orbit space can be given a structure of a classical orbifold, \cite{Sat56}, \cite{MOE}, \cite{MM03}, hence the name orbifold groupoid. However, the groupoid theory is more flexible than the orbifold theory. In particular, one often finds a global action groupoid in a Morita equivalence class of an orbifold. This is useful especially if the theory is applied for the noncommutative algebras where the localization techniques are complicated. The algebraic model of the quantum weighted projective spaces is constructed using this approach \cite{BF12}.

Associated with an orbifold there are two relevant complex function algebras: the algebra of smooth invariant functions and the smooth convolution algebra. The invariant subalgebra can be considered as a model describing properties of the quotient whereas the convolution algebra would describe the equivariant properties. For example, the geodesic length between a pair of points in the same orbit is equal to zero \cite{GH06}. This is an invariant property and one should be able to recover the metric properties from the invariant spectral triple, recall the geodesic length theorem \cite{Con94}. The well known shortcoming of this theory is that the universal differential calculus of an invariant spectral triple is not sufficiently rich. Namely, if the groupoid action on its base is not free, then the space of invariant differential forms cannot be reconstructed from the invariant spectral triple \cite{RV08}. The convolution algebras give rise to a proper differential calculus, homology and K-theory. For example this is demonstrated in the study of equivariant Fredholm index problems in \cite{GL13}. In this work we develop spectral triples over both algebras under certain natural conditions. 

\3 \noindent \textbf{I.2.} An orbifold is modelled as a Morita equivalence class of a groupoid and therefore one should have a precise understanding of how Morita equivalences operate on spectral triples. We approach this problem by fixing a proper {\' e}tale groupoid $X_{\bullet}$ as a representative of an orbifold. In addition, we assume that $X_{\bullet}$ is a spin groupoid since we are dealing with Dirac spectral triples. There is a Dirac spectral triple over the algebra of invariant functions $C^{\infty}(X)^{\Theta}$ if $X_{\bullet}$ is compact and there is a Dirac spectral triple over the groupoid convolution algebra $C^{\infty}_c(\Xi)$ if $X_{\bullet}$ is effective. A detailed description of these spectral triples is given in \cite{Har14}. 

Now the primary problem is to associate morphisms of spectral triples with the geometric Morita equivalences. If $X_{\bullet}$ represents an orbifold and $\phi$ is a Morita equivalence between $X_{\bullet}$ and $Y_{\bullet}$ then, under the compactness and effectiveness hypothesis discussed above, there is the following diagram:  
\begin{center}
\begin{tikzpicture}
\path (0,2) node (a) {$(C^{\infty}_c(\Theta), \eth, L^2(F_{\Sigma}))$} (6,2) node (b) {$(C^{\infty}_c(\Xi), \phi_{\#} \eth, L^2(\phi_{\#} F_{\Sigma}))$}  (3,1) node (c) {$X_{\bullet} \stackrel{\phi}{\longleftrightarrow} Y_{\bullet}$} (0,0) node (d) {$(C^{\infty}(X)^{\Theta}, \eth,L^2(F_{\Sigma})^{\Theta})$}  (6,0) node (e) {$(C^{\infty}(Y)^{\Xi}, \phi_{\#} \eth, L^2(\phi_{\#} F_{\Sigma})^{\Xi})$}; 
\begin{scope}
\draw[densely dashed,->] (c) --  (a) ; 
\draw[densely dashed,->] (c) --  (b); 
\draw[densely dashed,->] (c) --  (d) ; 
\draw[densely dashed,->] (c) --  (e);
\draw[->] (a) -- node [above] {$\phi_{\#}$} (b) ; 
\draw[->] (d) -- node [below] {$U_{\phi}$} (e);  
\end{scope}
\end{tikzpicture}
\end{center}
The symbols $F_{\Sigma}$ and $\eth$ denote a spinor bundle on $X_{\bullet}$ and a $\Theta$-invariant Dirac operator. The dashed lines indicates that there is a spectral triple associated with the groupoid. On the level of invariant spectral triples $\phi$ determines a unitary equivalence of spectral triples $U_{\phi}$. This is a spectral geometric analogue of the property that the orbit space is Morita invariant. On the level of convolution algebras we get an induced spectral triple over $C^{\infty}_c(\Xi)$. Whenever the dimension of the groupoid base is an even integer the horizontal morphisms are compatible with the even structure determined by the chirality operator.  The construction is independent on the choice of the spin proper {\' e}tale groupoid $X_{\bullet}$ that represents the orbifold: if $\phi$ is an {\' e}tale structure preserving Morita equivalence, then the spinor bundles and Dirac operators on $Y_{\bullet}$ are, up to an equivalence, the induced spinor bundles and Dirac operators through $\phi$. 

A major open problem in noncommutative geometry is to undestand the functoriality of spectral triples.  In the case of Morita equivalences, a particularly challenging issue is to properly axiomatize the transformation of Dirac operators. An important property in the diagram is that the induced Dirac operator $\phi_{\#} \eth$ is $\Xi$-invariant and a Morita equivalence restricts to define a unitary equivalence of Dirac operators on the Hilbert spaces of invariant spinors. So, the induced Dirac operator is subject to very specific conditions which can be easily studied in the context of noncommutative spectral triples as well.  Any generalized homomorphism of Lie groupoids can be represented as a composition of a Morita equivalence and a strict groupoid homomorphism and so this theory can can be also adapted to the study of generalized homomorphisms of spectral triples.

\3 \noindent \textbf{I.3.} Consider the case of a compact action groupoid $G \times X \rightrightarrows X$ representing an orbifold. If the action of the group $G$ on $X$ is free the orbit space $X/G$ is a smooth manifold and there is a Morita equivalence between the action groupoid and the unit groupoid over $X/G$. The convolution algebra of the unit groupoid $X/G$ is the algebra $C^{\infty}(X/G)$ with the pointwise product. It thus follows that on the level of convolution spectral triples the diagram of I.2 restricts to
\begin{eqnarray*}
(C^{\infty}_c(G \times X), \eth, L^2(F_{\Sigma})) \stackrel{\phi_{\#}}{\longrightarrow} (C^{\infty}(X / G), \phi_{\#} \eth, L^2(\phi_{\#} F_{\Sigma}))
\end{eqnarray*}
so that $\phi_{\#} \eth$ and $L^2(\phi_{\#} F_{\Sigma})$ identify with the Dirac operator and the Hilbert space of spinors on the manifold $X/G$. Now the spectral triple on the right hand side is unitarily equivalent to the invariant spectral triple $(C^{\infty}(X)^G, \eth, L^2(F_{\Sigma})^G)$, which is a consequence of the lower horizontal arrow in the diagram of I.2. This suggest that freeness of an action in noncommutative geometry should be axiomatized as an existence of a Morita equivalence between the convolution algebra spectral triple and the invariant spectral triple. So, ultimately, in order to understand freeness in noncommutative geometry one needs to axiomatize the Morita equivalence of spectral triples over orbifolds.

\3 \noindent \textbf{I.4  Notation.} We use the symbol $X_{\bullet}$ to denote a Lie groupoid with a base manifold $X_{(0)}$ and a manifold of $k$-times composable arrows $X_{(k)}$ for all $k \in \mathbb{N}$. We also write occasionally $X_{(1)} = \Theta$ and $X_{(0)} = X$ and then use the symbol $\Theta \rightrightarrows X$ to denote the groupoid. $\Theta$ and $X$ are smooth manifolds which are both Hausdorff and second countable. The target and source maps: $t,s : \Theta \rightarrow X$ are smooth submersions. $X_{\bullet}$ has a simplicial structure. $\textbf{1}_{x}$ denotes the unit morphism at $x \in X$ on the base. If $U$ and $V$ are subsets in $X$ we define the following subspaces in $\Theta$: 
\begin{eqnarray*}
\Theta_U = s^{-1}(U), \5 \Theta^V = t^{-1}(V) \5 \text{and}\5 \Theta_U^V = s^{-1}(U) \cap t^{-1}(V)
\end{eqnarray*}
If $x \in X$ and $U = V = \{x\}$ then $\Theta_x^x$ is the isotropy group at $x$. The orbit space (coarse moduli space) of $X_{\bullet}$ is denoted by $|\Theta|$. A groupoid is called compact if $|\Theta|$ is compact.

All groupoids are Lie groupoids in this work. The parameter $n$ will denote the dimension of the base manifold $X$ everywhere below. We will frequently consider a case where $X_{\bullet}$ is a groupoid and $Y_{\bullet}$ is another groupoid which is Morita equivalent to $X_{\bullet}$ (to be defined below). Then we call $Y_{\bullet}$ a representative of $X_{\bullet}$. In addition, without an exception, we use the symbol $\Theta$ to denote the manifold of arrows in $X_{\bullet}$ and $\Xi$ to denote the manifold of arrows in $Y_{\bullet}$.  \3

This project is funded by the John Templeton Foundation and "SFB 1085 Higher Invariants". 

\section{Morita Equivalence of Groupoids}

\noindent \textbf{1.1.} Let $X_{\bullet}$ and $Y_{\bullet}$ denote a pair of Lie groupoids. A generalized homomorphism between $X_{\bullet}$ and $Y_{\bullet}$ is a triple $(\varrho, Q,  \alpha)$ so that $Q$ is a smooth manifold, $\varrho$ and $\alpha$ are smooth maps in the diagram
\begin{center}
\gpdhom{Q}{\Theta}{X}{Y}{\Xi}{\varrho}{\alpha}
\end{center}
and, in addition:
\begin{quote}
\textbf{1.} $\varrho$ is an anchor for a left action of $\Theta$ on $Q$, $\alpha$ is an anchor for a right action of $\Xi$ on $Q$, and the actions are mutually commutative.

\textbf{2.} $\varrho$ is a surjective submersion and the action of $\Xi$ on $Q$ is free and transitive on the fibres of $\varrho$: if $q,q' \in Q$ and $\varrho(q) = \varrho(q')$ then there exists a unique $\tau \in \Xi$ so that $q \cdot \tau = q'$. 
\end{quote}

\noindent The condition 2 means that the right action of $\Xi$ on $Q$ makes $Q$ a right $\Xi$-torsor over $X$. If $(\varrho_1, Q_1, \alpha_1)$ and $(\varrho_2, Q_2, \alpha_2)$ are generalized homomorphisms $X_{\bullet} \leftrightarrow Y_{\bullet}$ and $Y_{\bullet} \leftrightarrow Z_{\bullet}$, then they can be composed to a generalized homomorphism  $X_{\bullet} \leftrightarrow Z_{\bullet}$ by setting
\begin{eqnarray*}
\big[ Q_1 \fp{\alpha_1}{\varrho_2} Q_2 \big] / \sim \5 \text{with} \5 (q_1,q_2) \sim (q_1 \cdot \sigma, \sigma^{-1} \cdot q_2)
\end{eqnarray*}
(the symbol $\fp{\alpha_1}{\varrho_2}$ denotes a fibre product). The anchor maps for the composition are $\varrho_1$ for the left action and $\alpha_2$ for the right action. The composition is associative. The following is useful \cite{LTX07}.\3
\noindent \textbf{Proposition 1.} Let $Y_{\bullet}$ be a Lie groupoid and $X$ a manifold. Then the following are equivalent: 
\begin{quote}
\textbf{1.} $Q$ is a right $\Xi$-torsor over $X$. 

\textbf{2.} $Q$ is a right $\Xi$-space where the $\Xi$-action is free and proper so that $X \simeq Q/\Xi$.
\end{quote}

A 2-morphism, or an equivalence, between two generalized homomorphism $(\varrho_1, Q_1, \alpha_1)$ and $(\varrho_2, Q_2, \alpha_2)$ is a $\Theta-\Xi$ equivariant diffeomorphism $T: Q_1 \rightarrow Q_2$. Lie groupoids, generalized homomorphisms and 2-morphisms define a 2-category. 

A Morita equivalence  is a generalized homomorphism $\phi = (\varrho, Q, \alpha)$ so that both actions make $Q$ a groupoid torsor: in addition to the properties above, the anchor $\alpha$ is a surjective submersion and the $\Theta$-action on $Q$ makes $Q$ a left $\Theta$-torsor over $Y$: if $q,q' \in Q$ and $\alpha(q) = \alpha(q')$, then there exists a unique $\sigma \in \Theta$ so that $\sigma \cdot q = q'$. In the case of a Morita equivalence, $Q$ is called a $X_{\bullet}-Y_{\bullet}$ bitorsor. Morita equivalences are not unique but two Morita bitorsors between the same groupoids can be mapped to each other under 2-morphisms. \3

\noindent \textbf{1.2.} A groupoid $\Theta$ is defined to be an {\' e}tale groupoid if its source and target maps are both local diffeomorphisms. It is a proper groupoid if $(s,t): \Theta \rightarrow X \times X$ is a proper map. The property of being a proper groupoid is invariant under the Morita equivalences, however, being {\' e}tale is not invariant.  A groupoid which is Morita equivalent to a proper {\' e}tale groupoid is called an orbifold groupoid. \3

\noindent \textbf{Proposition 2.} Suppose that $X_{\bullet}$ is an {\' e}tale Lie groupoid and $(\varrho,Q,\alpha)$ is a Morita equivalence between $X_{\bullet}$ and $Y_{\bullet}$, then $\alpha: Q \rightarrow Y$ is a local diffeomorphism. If the bitorsor $Q$ is compact, then  
\begin{quote}
\textbf{1.} $\alpha: Q \rightarrow Y$ is a smooth covering projection. 

\textbf{2.} $\varrho: Q \rightarrow X$ is a smooth fibre bundle projection.
\end{quote}
\noindent Proof. Since $X_{\bullet}$ is {\' e}tale, the source map $s: \Theta \rightarrow X$ has discrete fibres. If $y \in Y$, then the fibre $\alpha^{-1}(y)$ can be identified with $\Theta_x$ as a set if $x \in X$ is such that $\varrho(q) = x$ for some  $q \in \alpha^{-1}(y)$. It follows that the dimension of the fibres of $\alpha: Q \rightarrow Y$ are equal to the dimension of $Y$. Since $\alpha: Q \rightarrow Y$ is a submersion, it is also an immersion. It follows that $\alpha$ is a local diffeomorphism. 
If $Q$ is also compact, then $\alpha$ and $\varrho$ are proper maps. Then 1 and 2 follow from Ehresmann's theorem. \5 $\square$\3

By the symmetry of the Morita equivalence one observes that if $X_{\bullet}$ and $Y_{\bullet}$ are both proper and {\' e}tale, then the maps $\alpha$ and $\varrho$ are both local diffeomorphisms. However, if $Y_{\bullet}$ is not {\' e}tale, then $\varrho$ does not need to be a local diffeomorphism. In particular, the dimension of a groupoid base is not Morita invariant.\3
 
\noindent \textbf{1.3.} A vector bundle on a Lie groupoid $X_{\bullet}$ is a smooth vector bundle on the groupoid base, $\pi: \xi \rightarrow X$ with a typical fibre $V$ subject to an action of the groupoid from the left: 
\begin{eqnarray*}
\Theta \fp{s}{\pi} \xi \rightarrow \xi, \5 (\sigma, u_x) \mapsto (\rho(\sigma)u)_{\sigma \cdot x}
\end{eqnarray*}
where $\rho: \Theta \rightarrow GL(V)$ satisfies 
\begin{eqnarray*}
\rho(\tau) \rho(\sigma) = \rho(\tau \sigma), \5 \rho(\textbf{1}_x) = \iota
\end{eqnarray*}
for all $(\tau, \sigma) \in X_{(2)}$ and unit arrows $\textbf{1}_x$. The vector space $V$ can be real or complex but we always assume that it has a finite rank. Two vector bundles on $X_{\bullet}$ are defined to be isomorphic if they are isomorphic as geometric vector bundles on $X$ and the vector bundle isomorphism $\theta: \xi_1 \rightarrow \xi_2$ commutes with the actions of $\Theta$ on $\xi_1$ and $\xi_2$. 

The Morita equivalences can be used to transport vector bundles between Lie groupoids. The similar strategy was used in \cite{LTX07} for principal bundles. If $\xi$ is a vector bundle on $X_{\bullet}$ and $\phi = (\varrho, Q,\alpha)$ is a Morita equivalence $X_{\bullet} \leftrightarrow Y_{\bullet}$, then the induced bundle $\phi_{\#} \xi$ is the vector bundle on $Y_{\bullet}$ defined by 
\begin{eqnarray*}
\phi_{\#} \xi = \big[Q \fp{\varrho}{\pi} \xi \big] / \Theta
\end{eqnarray*}
where the $\Theta$-action is given by 
\begin{eqnarray*}
\Theta \fp{s}{\varrho} \phi_{\#} \xi \rightarrow \phi_{\#} \xi, \5 \sigma \cdot [q,u] = [\sigma \cdot q, \rho(\sigma)u].   
\end{eqnarray*}
The bundle projection sends $[q,u]$ to $\alpha(q)$. The fibre of $\phi_{\#} \xi$ at $y \in Y$ is an equivalence class of fibres of $\xi$ parametrized by the points $\alpha^{-1}(y)$ in $Q$, and if $q \in \alpha^{-1}(y)$, then the fibre is presented by the fibre of $\xi$ at $\varrho(q)$.  The action of $\Xi$ on $\phi_{\#} \xi$ is defined by 
\begin{eqnarray*}
[q,u] \cdot \tau \mapsto [q \cdot \tau , u]
\end{eqnarray*} 
for all $\tau$ such that $t(\tau) = \alpha(q)$.

Consider $X_{\bullet}$ and $Y_{\bullet}$ to be Morita equivalent and $\phi^1 = (\varrho_1, Q_1, \alpha_1)$ and $\phi^2 = (\varrho_2, Q_2, \alpha_2)$ a pair of Morita equivalences with a $\Theta - \Xi$ equivariant diffeomorphism $T: Q_1 \rightarrow Q_2$. Then the induced bundles $\phi^1_{\#} \xi$ and $\phi^2_{\#} \xi$ are equivalent and there is a canonical vector bundle equivalence defined by 
\begin{eqnarray*}
\phi^1_{\#} \xi \rightarrow \phi^2_{\#}\xi; \5 [q,u] \mapsto [T(q), u]
\end{eqnarray*}
where $q \in Q$ and $u \in \xi_{\varrho_1(q)} = \xi_{\varrho_2(T(q))}$.\3

\noindent \textbf{1.4.} Isomorphism classes of vector bundles on a proper {\' e}tale groupoid $X_{\bullet}$ can be described in terms of cohomology theoretic data, \cite{Har14}. We return to this in 1.5. In general, the groupoid cohomology is a cohomology of a double complex where one direction is determined by the complex associated to the simplicial structure and the other direction is determined by an injective resolution. We shall only employ cohomology in degree one and for this reason we will introduce a localization procedure which makes the injective resolution redundant when the degree one cohomology is considered, \cite{Har14}.  

Let $X_{\bullet}$ be a Lie groupoid. Fix an open cover $\{N_a\}$ for the groupoid base $X$. Associated to the cover there is a {\v C}ech groupoid  $\check{X}_{\bullet}$ defined as follows: the arrows and the base of $\check{X}_{\bullet}$ are given by 
\begin{eqnarray*}
\coprod_{ab} \Theta_{N_b}^{N_a} \rightrightarrows \coprod_a N_a 
\end{eqnarray*}
and the source and target maps are defined by restriction of the source an target maps on $X_{\bullet}$. The {\v C}ech groupoid $\check{X}_{\bullet}$ is Morita equivalent to $X_{\bullet}$.

Consider a Morita equivalence $\phi = (\varrho, Q, \alpha)$ between $X_{\bullet}$ and $Y_{\bullet}$. In what follows we fix open covers for the base manifolds $X$ and $Y$  and study how to localize the Morita equivalence $\phi$ accordingly. This allows us to transport locally defined data from $X_{\bullet}$ to $Y_{\bullet}$. Since $X$ and $Y$ are manifolds we can equip them with good covers $\{U_a\}$ and $\{V_i\}$ respectively. We do not write the index sets explicitly, but in both cases the index set can be chosen to be countable since we are working with manifolds. Then define the bitorsor
\begin{eqnarray*}
\check{Q} = \coprod_{a,i} \varrho^{-1}(U_{a}) \cap \alpha^{-1}(V_i).
\end{eqnarray*}
We let $\check{X}_{\bullet}$ and $\check{Y}_{\bullet}$ be the {\v C}ech groupoids
\begin{eqnarray*}
\coprod_{ab} \Theta_{ab} \rightrightarrows \coprod_a U_a, \5 \text{and} \5 \coprod_{ij} \Xi_{ij} \rightrightarrows \coprod_i V_i
\end{eqnarray*} 
where we have written $\Theta_{ab} = \Theta_{U_b}^{U_a}$ and $\Xi_{ij} = \Xi_{V_j}^{V_i}$. Define the maps $\check{\varrho}: \check{Q} \rightarrow \check{X}$ and $\check{\alpha}: \check{Q} \rightarrow \check{Y}$ so that $\check{\varrho}$ sends $q \in \varrho^{-1}(U_{a}) \cap \alpha^{-1}(V_i)$ to $x \in U_{a}$ if $\varrho(q) = x$ and similarly $\check{\alpha}$ sends $q \in \varrho^{-1}(U_{a}) \cap \alpha^{-1}(V_i)$ to $y \in V_{i}$ if $\alpha(q) = y$.  \3

\noindent \textbf{Lemma 1.} In the notation of the previous paragraph, a Morita equivalence $\phi$ restricts to a Morita equivalence between the {\v C}ech groupoids $\check{X}_{\bullet} \leftrightarrow \check{Y}_{\bullet}$ which is given by
\begin{eqnarray*}
\check{\phi} = (\check{\varrho}, \check{Q}, \check{\alpha}).
\end{eqnarray*}
 
\noindent Proof. The maps $\check{\varrho}$ and $\check{\alpha}$ are surjective because $\varrho$ and $\alpha$ are. They are submersions since locally, in a small enough neighborhood, they coincide with $\varrho$ and $\alpha$. The map $\check{\varrho}$ is the anchor for the action of $\check{\Theta}$ and so the domain of the action is given by
\begin{eqnarray*}
[\coprod_{ab} \Theta_{ab}] \fp{s}{\check{\varrho}} [\coprod_{c,i} \varrho^{-1}(U_{c}) \cap \alpha^{-1}(V_i)] = \coprod_{ab,i} (\Theta_{ab} ) \fp{s}{\check{\varrho}} \varrho^{-1}(U_{b}) \cap \alpha^{-1}(V_i). 
\end{eqnarray*}
We shall write $q_{a,i}$ to indicate that $q \in \varrho^{-1}(U_{a}) \cap \alpha^{-1}(V_i)$. If $\sigma \in \Theta_{ab}$, then the left action is given in its domain by
\begin{eqnarray*}
(\sigma, q_{b,i}) \mapsto (\sigma \cdot q)_{a,i}. 
\end{eqnarray*}
This is well defined since the $\Theta$-action is in the direction of the fibres of $\check{\alpha}$ and therefore $\sigma \cdot q \in \alpha^{-1}(V_i)$. The groupoid $\check{\Theta}$ acts transitively and freely on the $\check{\alpha}$-fibres because $\Theta$ acts transitively and freely on the $\alpha: Q \rightarrow Y$ fibres. 

The right action of $\check{\Xi}$ has the anchor $\check{\alpha}$ and the domain of the action is given by 
\begin{eqnarray*}
 [\coprod_{a,i} \varrho^{-1}(U_{a}) \cap \alpha^{-1}(V_i)] \fp{\check{\alpha}}{t} [\coprod_{jk} \Xi_{jk}]  = \coprod_{a,ij} \varrho^{-1}(U_{a}) \cap \alpha^{-1}(V_i)  \fp{\check{\alpha}}{t} \Xi_{ij}. 
\end{eqnarray*}
If $\tau \in \Xi_{ij}$, then the right action is defined in its domain by
\begin{eqnarray*}
(q_{a,i}, \tau) \mapsto (q \cdot \tau)_{a,j}.
\end{eqnarray*}
The $\check{\Xi}$-action is transitive and free on the fibres because this is the case for the $\Xi$-action on $Q$. The left and the right actions are commutative because the corresponding actions of $\Theta$ and $\Xi$ on $Q$ are commutative. \5 $\square$\3

\noindent \textbf{1.5.} Let $X_{\bullet}$ be a proper {\' e}tale groupoid. We take $GL_k$ to be the group of nonsingular real or complex $k \times k$ matrices. Denote by $\tx{GL}_k$ the sheaf of smooth $GL_k$ valued maps. $GL_k$ is not an abelian group, however, the first degree sheaf cohomology groups valued in $\tx{GL}_k$ are well defined as usual, \cite{BX11}. Let $\{U_a\}$ be a good open cover of the base $X$ and let $\check{X}_{\bullet}$ denote the {\v C}ech groupoid. Associated with the simplicial structure of the groupoid $\check{X}_{\bullet}$ defined in I.4., there is the simplicial sheaf cohomology group $H^1(\check{X}_{\bullet}, \tx{GL}_k)$.  The isomorphism classes of (real or complex) vector bundles on $X_{\bullet}$ correspond to the classes of $H^1(\check{X}_{\bullet}, \tx{GL}_k)$.  This correspondence was described explicitly in \cite{Har14}:  each vector bundle on $X_{\bullet}$ defines a structure cocycle which is represented in $H^1(\check{X}_{\bullet}, \tx{GL}_k)$. The structure cocycle $g$ has the components
\begin{eqnarray*}
g : \prod_{ab} \Theta_{U_b}^{U_a} \rightarrow GL_k. 
\end{eqnarray*}
The cocycle condition reads $g(\tau)g(\sigma) = g(\tau \sigma)$ for all $(\tau, \sigma) \in \check{X}_{(2)}$. The ordinary geometric structure of the reconstructed bundle is determined by the values of the cocycle on the unit arrows $\textbf{1}_x$: let $U_a \cap U_b \neq \emptyset$, then $g(\textbf{1}_x) : \Theta_{U_b}^{U_a} \rightarrow GL_k$ defines the geometric transition functions which are applied to glue all the local components $U_a \times E$ and $U_b \times E$ over $U_a \cap U_b$ ($E$ is $\mathbb{R}^k$ of $\mathbb{C}^k$). The nonidentity arrows define the groupoid action: if $\sigma \in \Theta_{U_b}^{U_a}$, then the action is given locally, in the trivialization of $U_b$ by 
\begin{eqnarray*}
\sigma \cdot [x_b, u] = [(\sigma \cdot x)_a, g(\sigma) u] 
\end{eqnarray*}
where $s(\sigma) = x_b$. The right side is written in the trivialization on $U_a$. The unit arrows act now as identities since the geometric reconstruction has identified the vectors related by such actions. 

Suppose that $X_{\bullet}$ and $Y_{\bullet}$ are as in 1.4. In the following we need a section for the local diffeomorphism $\check{\alpha}: \check{Q} \rightarrow \check{Y}$ of 1.4. Each cover sheet $V_i$ can be chosen to be small enough so that $\alpha$ has a section on $V_i$. Moreover, by shrinking $V_i$ further if necessary, we can assume that the image of the section lies in one of the cover sheets $\varrho^{-1}(U_a) \cap \alpha^{-1}(V_i)$. Such a section $\beta$ has the local components over each $i$ given by
\begin{eqnarray}\label{beta}
\beta_a^i : V_i \rightarrow \varrho^{-1}(U_a) \cap \alpha^{-1}(V_i)
\end{eqnarray}
where $U_a$ is one of the cover sheets of $X$. \3

\noindent \textbf{Proposition 3.} Let $\check{X}_{\bullet}$ and $\check{Y}_{\bullet}$ be the {\v C}ech groupoids and $\check{\phi}$ the Morita equivalence of Lemma 1. Suppose that $\beta$ is the section of $\check{\alpha}$ with components \eqref{beta}. If $g \in H^1(\check{X}_{\bullet}, \tx{GL}_k)$, then
\begin{quote}
\textbf{1.}  A representative of the induced cocycle $\check{\phi}_{\#}g$ in $H^1(\check{Y}_{\bullet}, \tx{GL}_k)$ is given by 
\begin{eqnarray*}
(\check{\phi}_{\#}g)_{ij}(y_{j} \rightarrow y'_{i}) = g_{ab}(\varrho \circ \beta_b^j(y_{j}) \rightarrow \varrho \circ \beta_a^i(y'_{i}))
\end{eqnarray*}

\textbf{2.} If $g$ is a structure cocycle of the vector bundle $\xi$ on $X_{\bullet}$, then $\check{\phi}_{\#}g$ is a structure cocycle of the vector bundle $\phi_{\#} \xi$ on $Y_{\bullet}$. 
\end{quote}

\noindent Proof. Following the discussion in the appendix A.1, there is a groupoid $\check{Q}_{(1)} \rightrightarrows \check{Q}$ such that the Morita equivalence of Lemma 1  can be represented as a pair of weak equivalences 
\begin{eqnarray}\label{weakeq}
\check{X}_{\bullet} \leftarrow \check{Q}_{\bullet} \rightarrow \check{Y}_{\bullet}. 
\end{eqnarray} 
Given a sheaf cohomology class $g \in H^1(\check{X}_{\bullet}, \tx{GL}_k)$, we can use the first weak equivalence in \eqref{weakeq} to pull the class to the sheaf cohomology of $\check{Q}_{\bullet}$. Then we can use the section $\check{Y} \rightarrow \check{Q}$ to pull the class to the cohomology of $\check{Y}_{\bullet}$. This class is a representative of the induced cocycle in the cohomology of $\check{Y}_{\bullet}$ and the cohomology class is independent on the choice of the section, \cite{Beh04}. The first pullback gives a cocycle in $H^1(\check{Q}_{\bullet}, \tx{GL}_k)$ with the components: 
\begin{eqnarray*}
(\varrho^*g)_{ab}^{ij}(\gamma: q \rightarrow q') = g_{ab}(\varrho(q) \rightarrow \varrho(q')). 
\end{eqnarray*}
where $\gamma$ is an arrow $\varrho^{-1}(U_b) \cap \alpha^{-1}(V_j) \rightarrow \varrho^{-1}(U_a) \cap \alpha^{-1}(V_i)$. The pullback of this class through the section $\beta$ gives the cocycle $\check{\phi}_{\#}g$. 

For 2, we have the locally defined section $\beta_b^j$ defined over each component of $\{V_i\}$ and we use these to fix a representative for each fibre of $\phi_{\#} \xi$ over the components $V_{i}$. So, the fibre at $y_j \in V_j$ is represented by the fibre of $\xi$ at $\varrho \circ \beta_b^j(y_j)$ over $U_b$. Let $\Xi_{ij} \ni \tau : y_{j} \rightarrow y'_{i}$. The action by $\tau$ sends the vectors in the fibre at $y$ to vectors in the fibre at $y'$. Suppose that $\beta_b^j$ and $\beta_a^i$ are the components of $\beta$ over $V_j$ and $V_i$. Then the $\Xi$ action reads
\begin{eqnarray*}
\tau \cdot v(y_{j}) = g_{ab}(\varrho \circ \beta_b^j(y_{j}) \rightarrow \varrho \circ \beta_a^i(y'_{i})) v(y'_{i}),
\end{eqnarray*}
These are exactly the components of the pullback cocycle. \5 $\square$\3

In particular, given isomorphic bundles $\xi_1$ and $\xi_2$ on the groupoid $X_{\bullet}$, then the induced bundles $\phi_{\#} \xi_1$ and $\phi_{\#} \xi_2$ are isomorphic since bundle reconstruction depends only on the cohomology class of the transition data. \3

\noindent \textbf{1.6.} Let $X_{\bullet}$ and $Y_{\bullet}$ be as in 1.5 and let $\phi^1 = (\varrho_1, Q_1, \alpha_1)$ and $\phi^2 = (\varrho_2, Q_2, \alpha_2)$ be a pair of Morita equivalences $X_{\bullet} \leftrightarrow Y_{\bullet}$. It follows from the definition of the 2-morphisms that there are the equalities of sets: 
\begin{eqnarray*}
\varrho_1 \circ \alpha_1^{-1}(y) = \varrho_2 \circ \alpha_2^{-1}(y) \5 \text{and} \5 \alpha_1 \circ \varrho_1^{-1}(x) = \alpha_2 \circ \varrho_2^{-1}(x)
\end{eqnarray*}
for all $y \in Y$ and for all $x \in X$. The {\v C}ech groupoids $\check{X}_{\bullet}$ and $\check{Y}_{\bullet}$ are constructed  as in 1.4. Denote by $\check{Q}_1$ and $\check{Q}_2$ the {\v C}ech bitorsors associated to the Morita equivalences $\phi^1$ and $\phi^2$. 

Consider a vector bundle $\xi$ on $X_{\bullet}$. Associated to the Morita equivalences, $\phi^1$ and $\phi^2$, there are the vector bundles $\phi^1_{\#} \xi$ and $\phi^2_{\#} \xi$ on $Y_{\bullet}$. Let $\beta_a^i: \check{Y} \rightarrow \check{Q}_1$ be components of  a section of $\check{\alpha}_1$. Then $T \circ \beta_a^i: \check{Y} \rightarrow \check{Q}_2$ is a section of $\check{\alpha}_2$ since $\alpha_2 = \alpha_1 \circ T^{-1}$. It then follows that the structure cocycle of $\phi^2_{\#} \xi$ has a representative given by 
\begin{eqnarray*}
(\check{\phi}_2^{\#}g)_{ij}(y_{j} \rightarrow y'_{i}) &=& g_{ab}(\varrho_2 \circ T \circ \beta_b^j(y_{j}) \rightarrow \varrho_2 \circ T \circ \beta_a^i(y'_{i})) \\
&=& g_{ab}(\varrho_1  \circ \beta_b^j(y_{j}) \rightarrow \varrho_1 \circ \beta_a^i(y'_{i})).
\end{eqnarray*}
So, the transition data associated to these two choices of Morita equivalence are identical.

\section{Invariant Sections}

\noindent \textbf{2.1.} For now we denote by $C^{\infty}(X)$ the algebra of smooth functions on $X$ getting values either in $\mathbb{R}$ or $\mathbb{C}$. If $X_{\bullet}$ is a Lie groupoid then the algebra of invariant functions $C^{\infty}(X)^{\Theta}$ is the subalgebra of $C^{\infty}(X)$ of functions such that $f(x) = f(\sigma \cdot x)$ holds for all $x$ and $\sigma \in \Theta_x$. Alternatively, the invariance of a function can be stated by $s^*(f) = t^*(f)$.

Let $\phi = (\varrho, Q, \alpha)$ be a Morita equivalence $X_{\bullet} \leftrightarrow Y_{\bullet}$. There is an isomorphism of (real or complex) vector spaces 
\begin{eqnarray*}
&&\phi_{\#}: C^{\infty}(X)^{\Theta} \rightarrow C^{\infty}(Y)^{\Xi}. \\
&&\phi_{\#}(f): y \mapsto f(\varrho(q)) 
\end{eqnarray*}
where $q \in Q$ is any element in $\alpha^{-1}(y)$. Recall that $\Theta$ acts transitively in the direction of the fibres of $\alpha$: this implies that if $q'$ is another point so that $\alpha^{-1}(y) = q' \in Q$, then $q' = \sigma \cdot q$ for some $\sigma$ and 
\begin{eqnarray*}
f(\varrho(\sigma \cdot q)) = f(\sigma \cdot \varrho(q)) = f(\varrho(q))
\end{eqnarray*}
 by the $\Theta$-invariance. So, the choice of $q$ is arbitrary. The $\Xi$-invariance follows by construction: $\Xi$ acts in the direction of the the fibres of $\varrho$ but the isomorphism composes through $\varrho$. The isomorphism $\phi_{\#}$ will be studied with more details in a more general context below. We denote by $\phi^{-1}_{\#}$ the inverse of $\phi_{\#}$ which is defined by 
\begin{eqnarray*}
&&\phi_{\#}^{-1}: C^{\infty}(Y)^{\Xi} \rightarrow C^{\infty}(X)^{\Theta}. \\
&&\phi_{\#}^{-1}(g): x \mapsto g(\alpha(q)) 
\end{eqnarray*}
where $q$ is any point in $\varrho^{-1}(x)$. \3

\noindent \textbf{Note 1.} A general principle in the study of invariant structures on groupoids (such as invariant functions, vector bundle sections and operators) is that a Morita equivalence $X_{\bullet} \leftrightarrow Y_{\bullet}$ induces an equivalence $\phi_{\#}$ from a $\Theta$-invariant structure on $X_{\bullet}$ to the corresponding $\Xi$-invariant structure on $Y_{\bullet}$ so that the $\Theta$-invariance implies that $\phi_{\#}$ is well defined. Then the $\Xi$-invariance of the induced structures over the target groupoid $Y_{\bullet}$ will follow by construction since we always factor through the anchor map $\varrho$. \3

\noindent \textbf{2.2.} Let $X_{\bullet}$ be a groupoid and $\psi \in \Gamma(\xi)$ a smooth section of a vector bundle $\xi$ on $X_{\bullet}$. Then $\psi$ is defined to be $\Theta$-invariant if
\begin{eqnarray*}
\psi(x) = \rho^{-1}(\sigma) \psi(\sigma \cdot x)
\end{eqnarray*}
for all $\sigma \in \Theta$ and $x = s(\sigma)$.  Denote by $\Gamma(\xi)^{\Theta}$ the vector space of $\Theta$-invariant sections.

Next we consider a pair of Morita equivalent groupoids, $X_{\bullet}$ and $Y_{\bullet}$, and fix a Morita equivalence $\phi = (\varrho, Q, \alpha)$. The fibres of the bundle $\phi_{\#} \xi$ are equivalence classes of vector spaces and so we use the symbol $\psi_q(y)$ to denote that the section $\psi$ in $\phi_{\#} \xi$ gets its value at $y$ in the vector space at $q \in \alpha^{-1}(y)$: this vector space is the fibre of $\xi$ at $\varrho(q)$. Then, by definition
\begin{eqnarray*}
\rho^{-1}(\sigma)\psi_{\sigma\cdot q} (y) = \psi_{q}(y)
\end{eqnarray*}
for all $\sigma \in \Theta_{\varrho(q)}$, where $\rho$ determines the $\Theta$-action on $\xi$. \3
 
\noindent \textbf{Proposition 4.} Let $\xi$ be a vector bundle on $X_{\bullet}$. The linear map $\phi_{\#}$ which sends a section $\psi$ of $\xi$ to the section of $\phi_{\#} \xi$ given by 
\begin{eqnarray*}
\phi_{\#} \psi : y \mapsto [q,\psi(\varrho(q))] : = [\psi_q(\varrho(q))], \5 \text{$q \in Q$ is so that $\alpha(q) = y$}
\end{eqnarray*}
is an equivalence of vector spaces: 
\begin{eqnarray*}
\phi_{\#}: \Gamma(\xi)^{\Theta} \stackrel{\simeq}{\longrightarrow} \Gamma(\phi_{\#} \xi)^{\Xi}.
\end{eqnarray*}
In addition, if $f \in C^{\infty}(X)^{\Theta}$, then $\phi_{\#} f \psi = \phi_{\#}(f) \phi_{\#} \psi$. \3

\noindent Proof. Let us check that $\phi_{\#} \psi$ is well defined element in $\Gamma(\phi_{\#}\xi)^{\Xi}$. Let $q, \sigma \cdot q \in \alpha^{-1}(y)$. Then 
\begin{eqnarray*}
\rho^{-1}(\sigma)(\phi_{\#} \psi)_{\sigma\cdot q} (y) &=& \rho^{-1}(\sigma) \psi(\varrho(\sigma \cdot q)) \\
&=&  \rho^{-1}(\sigma) \psi(\sigma \cdot \varrho(q)) \\
&=& \psi(\varrho(q)) \\ &=& (\phi_{\#} \psi)_{q}(y). 
\end{eqnarray*}
In the third equality we used the $\Theta$-invariance of $\psi$. The $\Xi$-invariance of $\phi_{\#} \psi$, i.e. 
\begin{eqnarray*}
\phi_{\#}(\psi)(y) = \phi_{\#}(\psi)(\tau \cdot y), \5 \text{for all $\tau \in \Xi_y$},
\end{eqnarray*}
 holds since $\varrho$ is constant along the $\Xi$-orbits in $Q$.  Let $y \in Y$, $V$ be a neighborhood of $y$ and $\beta: V \rightarrow Q$ be a smooth section of $\alpha$. The local section exists because $\alpha$ is a local diffeomorphism. Then we can write locally 
\begin{eqnarray*}
(\phi_{\#} \psi)_{\beta(y)}(y) = \psi(\varrho \circ \beta(y)).  
\end{eqnarray*}
In particular, $\phi_{\#} \psi$ is smooth at $y \in Y$. Since $y$ is arbitrary, $\phi_{\#} \psi$ is smooth on $Y$.

Let $\zeta \in \Gamma(\phi_{\#}\xi)^{\Xi}$. Under the $\Theta$-action there are the relations 
\begin{eqnarray*}
\rho^{-1}(\sigma) \zeta_{\sigma \cdot q}(y) = \zeta_q(y).
\end{eqnarray*}
We define a section in $\xi$ by 
\begin{eqnarray*} 
A_{\phi} \zeta : x \mapsto \zeta_q(\alpha(q))
\end{eqnarray*}
where $q \in Q$ is any point so that $\varrho(q) = x$. The definition is independent on the choice of such $q$: any point in the same $\varrho$-fibre can be written by $q' = q \cdot \tau$ for some $\tau \in \Xi$ and the $\Xi$-invariance of $\zeta$ implies that $\zeta_{q\cdot \tau}(\alpha(q \cdot \tau)) = \zeta_q(\alpha(q))$. If $\sigma \in \Theta_x$, then 
\begin{eqnarray*}
\rho^{-1}(\sigma) A_{\phi} \zeta(\sigma \cdot x) &=& 
\rho^{-1}(\sigma) \zeta_{\sigma \cdot q}(\alpha(\sigma \cdot q)) \\ 
&=& \zeta_q(\alpha (\sigma \cdot q)) \\
&=&  \zeta_q(\alpha(q)) \\ &=& 
 A_{\phi} \zeta(x)
\end{eqnarray*}
where $q \in Q$ is such that $\varrho(q) = x$. Since $\varrho$ is a submersion it accepts local sections. If $x \in X$, $U$ is a local neighborhood of $x$ and $\delta: U \rightarrow Q$ is a smooth section of $\varrho$, then   
\begin{eqnarray*}
A_{\phi} \zeta(x) = \zeta_{\delta(x)}(\alpha \circ \delta(x))
\end{eqnarray*}
which is smooth at $x$.  Therefore $A_{\phi}$ defines a map $\Gamma(\phi_{\#} \xi)^{\Xi} \rightarrow \Gamma(\xi)^{\Theta}$. It is an easy matter to check that the linear maps $\phi_{\#}$ and $A_{\phi}$ are inverses of each other and that the isomorphism $\phi_{\#}$ respects the module structure under the action of the smooth invariant functions. \5 $\square$ \3

We shall write $\phi_{\#}^{-1} = A_{\phi}$ in what follows. 

Consider a pair of Morita equivalences $\phi^i = (\varrho_i, Q_i, \alpha_i)$ for $i = 1,2$ and a $\Theta-\Xi$ equivariant diffeomorphism $T: Q_1 \rightarrow Q_2$. If $\psi \in \Gamma(\xi)^{\Theta}$, then we get two invariant sections, $\phi_{\#}^1 \psi$ and $\phi_{\#}^2 \psi$ in the bundles $\phi^1_{\#} \xi$ and $\phi^2_{\#} \xi$. However, these sections are essentially identical: at $y \in Y$ we have the fibres $\alpha_1^{-1}(y)$ and $\alpha_2^{-1}(y)$ and the fibres of $\phi^1_{\#} \xi$ and $\phi^2_{\#} \xi$ are equivalence classes of fibres of $\xi$ in $\varrho_1 \circ \alpha_1^{-1}(y)$ and $\varrho_2 \circ \alpha_2^{-1}$. However, these are the same subset of $X$. So, for all $q_2 \in \alpha_2^{-1}(y)$ we have $q_1 \in \alpha_1^{-1}(y)$ such that $T(q_1) = q_2$. Now
\begin{eqnarray*}
\psi(\varrho_1(q_1)) = \psi(\varrho_2(T(q_1))  =  \psi(\varrho_2(q_2)).
\end{eqnarray*}
It then follows from the equivariance of $T$ that
\begin{eqnarray*}
(\phi^1_{\#} \psi)_{\sigma \cdot q_1}(y) = (\phi^2_{\#} \psi)_{T(\sigma \cdot q_1)}(y) 
\end{eqnarray*}
which implies that at $y \in Y$, the two induced sections are merely two description of the same section. This holds for all $y \in Y$. \3

\noindent \textbf{2.3.} For each arrow $\sigma \in \Theta_x$ there is an open neighborhood $U$ of $x$ and a local section of $s$, $\hat{\sigma}: U \rightarrow \Theta$, such that $\hat{\sigma}(x) = \sigma$ and $t \circ \hat{\sigma}$ is an open embedding. These sections are called local bisections and they always exist on Lie groupoids, \cite{MM03}. In the case of an {\' e}tale Lie groupoid, any two such sections agree on their common domain. The assignment
\begin{eqnarray*}
\sigma \mapsto \varphi_{\sigma} = t \circ \hat{\sigma}: U\rightarrow \varphi_{\sigma}(U). 
\end{eqnarray*}
defines a local diffeomorphisms. Denote by $\Delta(\Theta)$ the set of germs of these diffeomorphisms. An {\' e}tale groupoid is defined to be effective if the assignment sending arrows to the elements in $\Delta(\Theta)$ is injective. 

Let $X_{\bullet}$ be a proper {\' e}tale groupoid. If $\sigma \in \Theta$, then we get a unique element $\varphi_{\sigma} \in \Delta(\Theta)$. Therefore  there is a canonical choice of a differential $(d\varphi_{\sigma})$ for each $\sigma$ which are invertible matrices.  This amounts to define a tangent bundle $\tau X_{\bullet}$ on $X_{\bullet}$. This bundle is an ordinary tangent bundle $\pi: \tau X \rightarrow X$ which is equipped with an action of the groupoid given by 
\begin{eqnarray*}
(\sigma, v_x) \mapsto (d \varphi_{\sigma})_x v_{\sigma \cdot x},  
\end{eqnarray*}
which is defined in $\Theta \fp{s}{\pi} \tau X$. Similarly, the exterior bundles $\wedge^k \tau X_{\bullet}^*$ are the geometric exterior bundles on $X$ together with the the groupoid action 
\begin{eqnarray*}
\sigma \cdot w(x, v_1, \ldots, v_k) = w(\sigma \cdot x, (d \varphi_{\sigma})_x v_1, \ldots, (d \varphi_{\sigma})_x v_k)
\end{eqnarray*}
if $v_i \in  \tau X$. We denote by $\mathfrak{X}(X_{\bullet})^{\Theta}$ the space of vector fields which are invariant under the right actions of all arrows: 
\begin{eqnarray*}
 V(x) = (d \varphi_{\sigma})^{-1} V(\sigma \cdot x) 
\end{eqnarray*}
for all $\sigma \in \Theta$ and $x = s(\sigma)$.  Similarly, $\Lambda^k(X_{\bullet})^{\Theta}$ is the space of $k$-forms invariant under the action of all arrows.

Suppose that $Y_{\bullet}$ is a representative of a proper {\' e}tale groupoid $X_{\bullet}$ and $\phi = (\varrho, Q, \alpha)$ is a Morita equivalence. The induced differential forms on $Y_{\bullet}$ are defined in terms of coordinate system of $X$. We can apply a localization strategy, similar to 1.4, to write the forms in the coordinates of $Y$: this is parallel to the standard procedure in differential geometry to write the pullback forms in the coordinate system of the target manifold. Choose the covers $\{U_a\}$ for $X$ and $\{V_i\}$ for $Y$ which are fine enough so that $\varrho$ and $\alpha$ have local sections over the cover sheets. Choose the local sections $\beta_i: V_{i} \rightarrow Q$ of $\alpha$. Locally, on $V_{i}$, we can write the induced forms in the coordinates of $V_{i}$: 
\begin{eqnarray}\label{2}
\phi_{\#}\Omega|_{V_{i}} = (\varrho \circ \beta_i)^*\Omega.
\end{eqnarray}
for all $\Omega \in \Lambda(X_{\bullet})^{\Theta}$. Here we are exploiting the notation: now we have applied the local sections $\beta_i$ to choose a representative in the equivalence class of vector spaces in each fibre of $\phi_{\#} \xi$. For the verification that this provides a well defined global form we shall need the following technical result. \3

\noindent \textbf{Lemma 2.} Let $X_{\bullet}$ be a proper {\' e}tale Lie groupoid and $Y_{\bullet}$  a representative of its Morita equivalence class. Fix a Morita equivalence $\phi = (\varrho, Q, \alpha)$. 
\begin{quote}
\textbf{1.} For all $\varphi_{\sigma} \in \Delta(\Theta)$ and $q \in \varrho^{-1}(s(\sigma))$ there is a lift of $\varphi_{\sigma}$: a local diffeomorphism in $Q$ satisfying
\begin{eqnarray*}
\varrho \circ \hat{\varphi}_{\sigma} = \varphi_{\sigma}  \circ \varrho
\end{eqnarray*}
in an open neighborhood of $q$ which arises from the action of $\Theta$ on $Q$. The germ of the lift is unique. 

\textbf{2.} For all $\varphi_{\tau} \in \Delta(\Xi)$ and $r \in \alpha^{-1}(s(\tau))$ there is a lift of $\varphi_{\tau}$: a local diffeomorphism in $Q$ satisfying
\begin{eqnarray*}
\alpha \circ (\hat{\varphi}_{\tau})^{-1} = \varphi_{\tau} \circ \alpha
\end{eqnarray*}
in an open neighborhood of $r$ which arises from the action of $\Xi$ on $Q$. The germ of the lift is unique. 
\end{quote}

\noindent Proof. Let $\varphi_{\sigma} \in \Delta(\Theta)$. We fix an open neighborhood $O$ for $q \in Q$ which is small enough so that $\varrho(O)$ lies in the domain of $\varphi_{\sigma}$. Then define a local diffeomorphism on $O$ so that if $q' \in O$, then 
\begin{eqnarray*}
\hat{\varphi}_{\sigma}(q') = \hat{\sigma}(\varrho(q')) \cdot q',
\end{eqnarray*}
where $\hat{\sigma}$ is the unique local bisection of $\sigma$. This is a local diffeomorphism and it is determined by the action of $\Theta$ on $Q$ and the germ of the lift is unique among the diffeomorphism that are determined by the $\Theta$ action because the bisection is unique. 

Let $\varphi_{\tau} \in \Delta(\Xi)$. Now $\varphi_{\tau} = t \circ \hat{\tau}$ for the local bisection $\hat{\tau}$. Then we can define the lifting as in the case 1 by applying the right action of $\Xi$ on $Q$. Now the uniqueness follows because $\alpha$ and $\varphi_{\tau}$ are local diffeomorphism, and so a solution of the lifting problem has necessary a unique germ. \5 $\square$ \3

\noindent \textbf{Proposition 5.} Let $X_{\bullet}$ be a proper {\' e}tale groupoid and $Y_{\bullet}$ a representative of $X_{\bullet}$. In the coordinate system of $Y_{\bullet}$, $\phi_{\#}$ is the linear map $\Lambda(X_{\bullet})^{\Theta} \rightarrow \Lambda(Y_{\bullet})^{\Xi}$ given by 
\begin{eqnarray*}
\phi_{\#} (f_0 df_1 \wedge \cdots \wedge df_k) = \phi_{\#}(f_0) d \phi_{\#}(f_1) \wedge \cdots \wedge d \phi_{\#}(f_k).
\end{eqnarray*}
The $\Xi$-invariance of the forms on $Y_{\bullet}$ means the invariance with respect to all local diffeomorphisms $\Delta(\Xi)$. \3

\noindent Proof. Using the realization \eqref{2} and the rules in calculus of differential forms we get a local equality
\begin{eqnarray*}
 (\varrho \circ \beta_i)^* (f_0df_1 \wedge \cdots \wedge df_k) =  f_0( \varrho \circ \beta_i) df_1(\varrho \circ \beta_i) \wedge \cdots \wedge df_k(\varrho \circ \beta_i) 
\end{eqnarray*}
which proves that  $\phi_{\#}$ is given in the coordinates of $Y$ is as claimed. 

Next we show that the forms written in the coordinates of $Y$ are well defined global forms. Consider an intersection $V_{i} \cap V_{j} \neq \emptyset$ and the local sections $\beta_i: V_{i} \rightarrow Q$ and $\beta_j: V_{j} \rightarrow Q$. The sections $\beta_i$ and $\beta_j$ do not need to agree on the intersection $V_i \cap V_j$. However, since $\Theta$ is transitive in the fibres of $\alpha$, the values of these sections need to be related by the action of $\Theta$ at each point in  $V_i \cap V_j$. It follows from the uniqueness of the lifts that there is $\hat{\varphi}_{\sigma}$ so that $\beta_i = \hat{\varphi}_{\sigma} \circ \beta_j$ holds in a neighborhood of $y \in V_i \cap V_j$. It then follows that 
\begin{eqnarray*}
(\varrho \circ \beta_i)^* \Omega|_y &=& (\varrho \circ \hat{\varphi}_{\sigma} \circ \beta_j)^* \Omega|_y \\
&=& (\varphi_{\sigma}  \circ \varrho  \circ \beta_j)^* \Omega|_y \\
&=& (\varrho  \circ \beta_j)^* \varphi_{\sigma}^* \Omega|_y \\
&=& (\varrho  \circ \beta_j)^* \Omega|_y 
\end{eqnarray*} 
the last equality follows from the $\Theta$-invariance of the form. So the form is well defined.

Let $\varphi_{\tau}$ be an element of $\Delta(\Xi)$ and $\tau: y \rightarrow \tau \cdot y$ be an arrow so that $y \in V_{j}$ and $\tau \cdot y \in V_{i}$. Let $\hat{\varphi}_{\tau}$ be a lift of $\varphi_{\tau}$ which sends the point $q = \beta_i(\tau \cdot y) \in Q$ to $q \cdot \tau$. If $\beta_j$ is a locally defined on $V_j$, it might happen that $\beta_i \circ \varphi_{\tau}(y) \neq (\hat{\varphi}_{\tau})^{-1} \circ \beta_j(y)$ but these points need to be in the same $\alpha$ fibre and so there is an arrow $\sigma$ so that its action on $Q$ sends $\beta_i \circ \varphi_{\tau}(y)$ to $(\hat{\varphi}_{\tau})^{-1} \circ \beta_j(y)$. Then we use the lift $\hat{\varphi}_{\sigma}$ to define 
 \begin{eqnarray*}
(\beta')_i = \hat{\varphi}_{\sigma} \circ \beta_i
\end{eqnarray*}
locally. This defines is a section of $\alpha$ in a small enough neighborhood of $\tau \cdot y$ since the local diffeomorphism $\hat{\varphi}_{\sigma}$ operates using the $\Theta$-action, i.e., in the direction of the fibres. By the first part of the proof, we can use any local section to pullback $\Theta$-invariant forms. So we will use $(\beta')_i$ in a neighborhood of $\tau \cdot y$ which gives 
\begin{eqnarray*}
(\beta')_i \circ \varphi_{\tau} = (\hat{\varphi}_{\tau})^{-1} \circ \beta_j.
\end{eqnarray*}
in a neighborhood of $y \in Y$. Then a straightforward computation gives
\begin{eqnarray*}
\varphi^*_{\tau} \circ (\varrho \circ (\beta')_i)^* &=& (\varrho \circ (\beta')_i \circ \varphi_{\tau})^* \\ 
&=&  (\varrho \circ  (\hat{\varphi}_{\tau})^{-1} \circ \beta_j )^* \\
&=& (\varrho\circ \beta_j )^*
\end{eqnarray*}
The last equality holds because the local diffeomorphism $(\hat{\varphi}_{\tau})^{-1}$ applies arrows in $\Xi$ and so operates in the direction of $\varrho$-fibres. Consequently, the pullback forms are invariant under pullbacks by local diffeomoprhisms which are determined by any choice of a local bisection. \5 $\square$ \3 

\noindent \textbf{Corollary 1.} The linear map $\phi_{\#}: \Lambda(X_{\bullet})^{\Theta} \rightarrow \phi_{\#} \Lambda(X_{\bullet})^{\Theta}$ is an isomorphism. \3

\noindent Proof. Consider the image of $\phi_{\#}$ in  $\Lambda(Y_{\bullet})^{\Xi}$. Recall that we have fixed an open cover $\{U_a\}$ for $X$. The inverse $\phi_{\#}^{-1}: \phi_{\#} \Lambda(X_{\bullet})^{\Theta} \rightarrow \Lambda(X_{\bullet})^{\Theta}$ can be defined by choosing a local section $\delta_a: U_a \rightarrow Q$ for $\varrho$ over each $U_a$ and then defining the local forms 
\begin{eqnarray*}
\phi_{\#}^{-1} \Phi|_{U_a} = (\alpha \circ \delta_a)^* \Phi.
\end{eqnarray*}
Let us check that this is a left inverse of $\phi_{\#}$. Suppose that $y \in V_i \subset Y$ and $\varrho \circ \beta_i(y) = x \in U_a \subset X$. Then
\begin{eqnarray*}
\phi_{\#}^{-1} \circ \phi_{\#}(\Omega)|_x = (\alpha \circ \delta_a)^* (\varrho \circ \beta_i)^*\Omega|_x = (\varrho \circ \beta_i \circ \alpha \circ \delta_a)^*\Omega|_x.
\end{eqnarray*}
Now we can write locally $\beta_i \circ \alpha = \hat{\varphi}_{\sigma}$ for some lift. Then 
\begin{eqnarray*}
\phi_{\#}^{-1} \circ \phi_{\#}(\Omega)|_x &=& (\varrho \circ \hat{\varphi}_{\sigma} \circ \delta_a)^*\Omega|_x \\
&=& (\varphi_{\sigma} \circ \varrho \circ \delta_a)^*\Omega|_x \\
&=& \varphi_{\sigma}^*(\Omega)|_x = \Omega_x. 
\end{eqnarray*}
This holds for all $x \in X$. The proof that $\phi_{\#}^{-1}$ is the right inverse of $\phi_{\#}$ is similar. \5 $\square$ \3

Consider the case of differential forms valued in some bundle $\xi$ on $X_{\bullet}$. Then we have the space of invariant sections $\Lambda^*(X_{\bullet}, \xi)^{\Theta}$. If $\phi$ is a Morita equivalence, then we can use Proposition 5 to identify the space $\phi_{\#} \Lambda^*(X_{\bullet}, \xi)^{\Theta}$ as a subspace in $\Lambda^*(Y_{\bullet}, \phi_{\#} \xi)^{\Xi}$ so that the $\Xi$ action in the latter is given by 
\begin{eqnarray*}
\varphi^*_{\tau} \Phi(y,v_1, \ldots, v_k) = \Phi(\varphi_{\tau}(y), (d \varphi_{\tau})v_1, \ldots, (d \varphi_{\tau})v_k) 
\end{eqnarray*}
where $\varphi_{\tau}$ is a local diffeomoprhism associated to $\tau$. From now on we always use this identification when induced forms differential forms are considered.

\section{Invariant Connections}

\noindent \textbf{3.1.} Let $X_{\bullet}$ be a proper {\' e}tale groupoid, $\xi$ a vector bundle on $X_{\bullet}$ and $Y_{\bullet}$ a representative of the Morita equivalence class of $X_{\bullet}$. A connection in the vector bundle $\xi$ is an ordinary geometric connection $\nabla: \Gamma(\xi) \rightarrow \Lambda^1(X, \xi)$ which is invariant under the locally defined pullback actions by the element of $\Delta(\Theta)$:  
\begin{eqnarray*}
\varphi^*_{\sigma} \nabla (\varphi^*_{\sigma})^{-1} = \nabla,
\end{eqnarray*}
holds locally, for all $\sigma \in \Theta$. A groupoid connection restricts to a linear map 
\begin{eqnarray*}
\nabla: \Gamma(\xi)^{\Theta} \rightarrow \Lambda^1(X_{\bullet}, \xi)^{\Theta}.  
\end{eqnarray*}
Fix a Morita equivalence $\phi = (\varrho, Q, \alpha)$ between $X_{\bullet}$ and $Y_{\bullet}$. The bundle $\phi_{\#} \xi$ can be equipped with the induced connection
\begin{eqnarray*}
\phi_{\#}\nabla = \phi_{\#} \circ \nabla \circ \phi_{\#}^{-1}. 
\end{eqnarray*}
Notice that since $\nabla$ is a map $\Gamma(\xi)^{\Theta} \rightarrow \Lambda^1(X_{\bullet}, \xi)^{\Theta}$, then $\phi_{\#}^{-1}$ on the right side of $\nabla$ is a map $\Gamma(\phi_{\#} \xi)^{\Xi} \rightarrow \Gamma(\xi)^{\Theta}$ and $\phi_{\#}$ on the left side of $\nabla$ is a map $\Lambda^1(X_{\bullet},\xi)^{\Theta} \rightarrow \Lambda^1(Y_{\bullet}, \phi_{\#}\xi)^{\Xi}$.\3

\noindent \textbf{Proposition 6.} The induced connection $\phi_{\#}\nabla: \Gamma(\xi)^{\Theta} \rightarrow \Lambda^1(Y_{\bullet}, \phi_{\#} \xi)^{\Xi}$ restricts to a groupoid connection on the invariant sections and satisfies 
\begin{eqnarray}\label{connection}
\phi_{\#}\nabla_{V}( \phi_{\#}\psi) = \phi_{\#}(\nabla_{\phi_{\#}^{-1} V}(\psi)). 
\end{eqnarray}
for all invariant vector fields $V \in \mathfrak{X}(Y)^{\Xi}$ and sections $\psi \in \Gamma(\xi)^{\Theta}$. \3 

\noindent Proof. Clearly $\phi_{\#}\nabla$ is linear. The Leibnitz rule holds: for all $f \in C^{\infty}(Y)^{\Xi}$ and $\psi \in \Gamma(\phi_{\#} \xi)^{\Xi}$
\begin{eqnarray*}
(\phi_{\#}\nabla)(f \psi) &=& \phi_{\#} \circ \nabla((\phi_{\#}^{-1}f)(\phi_{\#}^{-1}\psi)) \\
&=& \phi_{\#}((d \phi_{\#}^{-1}f) \otimes (\phi_{\#} \psi) + (\phi_{\#}^{-1} f) \nabla(\phi_{\#}^{-1}\psi)) \\
&=& df \otimes \psi + f (\phi_{\#} \nabla)\psi. 
\end{eqnarray*}
The last equality follows from Proposition 5. The $\Xi$-invariance holds by construction. Thus, $\phi_{\#} \nabla$ restricts to a groupoid connection on the invariant sections.   

Let $V \in \mathfrak{X}(Y)^{\Xi}$. Recall from 2.3 that we have a local coordinate expression for the induced forms. Let us fix local sections $\beta_i$ on $V_{i}$ as in 2.3. The fields $V$ and $\phi_{\#}^{-1} V$ are $\varrho \circ \beta_i$ related over $V_{i}$. It follows that 
\begin{eqnarray*}
i_V \circ \phi_{\#} = \phi_{\#} \circ i_{\phi^{-1}_{\#} V}
\end{eqnarray*}
and consequently, 
\begin{eqnarray*}
i_V \phi_{\#} \nabla &=& i_V(\phi_{\#} \circ \nabla \circ \phi_{\#}^{-1}) \\ &=& \phi_{\#} \circ i_{\phi_{\#}^{-1} V} \circ \nabla \circ \phi_{\#}^{-1} \\
&=& (\phi_{\#} \nabla)_{\phi_{\#}^{-1} V}. \5 \square 
\end{eqnarray*}
 
\noindent \textbf{3.2.} Suppose that $X_{\bullet}, Y_{\bullet}, \phi$ and $\xi$ are as in 3.1 and that $\xi$ is equipped with an $\Theta$-invariant inner product $(\cdot, \cdot)$. Recall that the fibre of $\phi_{\#} \xi$ at $y \in Y$ is an equivalence class of vector spaces $\xi_x$ for those $x \in X$ such that there is $q \in Q$ for which $\alpha(q) = y$ and $\varrho(q) = x$. The induced bundle $\phi_{\#} \xi$ can be equipped with the inner product given at $y \in Y$ by
\begin{eqnarray*}
([v_1], [v_2])_{\#,y} = (v^x_1, v^x_2)_{x}
\end{eqnarray*}
where $x \in X$ is any point which is related to $y$ as  in the discussion above, and $v_i^x$ is the vector in the equivalence class $[v_i]$ evaluated in the fibre at $x$. The choice of such $x$ is arbitrary since the inner product $(\cdot, \cdot)$ is $\Theta$-invariant and the equivalence class of the fibres of $\phi_{\#} \xi$ is parametrized by the orbits under the $\Theta$-action. The $\Xi$-invariance of the induced inner product holds since the $\Xi$-action is in the direction of $\varrho$-fibres in $Q$. 

The inner product $(\cdot, \cdot)_{\#}$ in $\phi_{\#}\xi$ can be applied to define a linear pairing of the invariant sections $\Gamma(\phi_{\#} \xi)^{\Xi} \otimes \Gamma(\phi_{\#} \xi)^{\Xi} \rightarrow C^{\infty}(Y)^{\Xi}$ which satisfies
\begin{eqnarray*}
(\phi_{\#} \psi_1, \phi_{\#} \psi_2)_{\#,y} = (\psi_1, \psi_2)_x = \phi_{\#}(\psi_1, \psi_2)_y. 
\end{eqnarray*}
where $x$ is related to $y$ as in the discussion above. The last equality holds since $(\psi_1, \psi_2)$ is a $\Theta$-invariant function on $X$. So we conclude that
\begin{eqnarray*}
(\phi_{\#} \psi_1, \phi_{\#} \psi_2)_{\#} = \phi_{\#} (\psi_1, \psi_2), 
\end{eqnarray*}
whenever $\phi_{\#} \psi_1, \phi_{\#} \psi_2 \in \Gamma(\phi_{\#} \xi)^{\Xi}$. Since $\phi_{\#}$ is an invertible map, it induces a one to one correspondence between linear pairings in $\Gamma(\xi)^{\Theta}$ and $\Gamma(\phi_{\#} \xi)^{\Xi}$.

\section{Spinors and Dirac Operators}

We first follow \cite{Har14} and recall the construction of spinor bundles and invariant Dirac operators on a proper {\' e}tale groupoid. Then in 4.2 the construction is extended to the Morita equivalence classes.\3
 
\noindent \textbf{4.1.} Let $X_{\bullet}$ be a proper {\' e}tale groupoid and let $n$ denote its dimension. The tangent bundle $\tau X_{\bullet}$ can be equipped with a $\Theta$-invariant riemannian structure so that the groupoid acts on the tangent spaces orthogonally. Let us choose a good cover for the base $X$ and define the associated {\v C}ech groupoid  $\check{X}_{\bullet}$, recall 1.4. Fix a structure cocycle $g$ for the tangent bundle $\tau X_{\bullet}$ which defines a class in $H^1(\check{X}_{\bullet}, \tx{O}_n)$. The groupoid is defined to be orientable if the structure group can be reduced to $SO_n$, i.e. $g$ defines an element in $H^1(\check{X}_{\bullet}, \tx{SO}_n)$.

Suppose that $X_{\bullet}$ is orientable. For each $x \in X$ there is a Clifford algebra associated with the $n$-dimensional tangent space and to the inner product determined by the riemannian structure. The group $\text{Spin}(n)$ lies in the Clifford algebras and the $\text{SO}_n$-valued transition cocycle $g$ of $\tau X_{\bullet}$ can be lifted to a $\text{Spin}_n$-valued cochain $\hat{g}$. Now $\text{Ad}(\hat{g})$ is a cocycle in $H^1(\check{X}_{\bullet}, \tx{SO}(\text{cl}_n))$ since the center of the Spin group vanishes under the adjoint action. Denote by $\text{CL}(X_{\bullet})$ the bundle of  complexified Clifford algebras reconstructed from this cocycle, recall 1.5. The groupoid $X_{\bullet}$ is defined to be spin if the second Stiefel-Whitney class associated with the lifted cochain $\hat{g}$ vanishes: this is the obstruction class for the lifting of the structure cocycle. If $X_{\bullet}$ is spin, then one can reconstruct a bundle of spinors from the cocycle $\hat{g}$: the fibres are spinor modules $(\rho_s, \Sigma)$, i.e. irreducible representations for the Clifford algebras with an action $\rho_s$ of $\text{Spin}_n$ on $\Sigma$ determined by the embedding of $\text{Spin}_n$ to the Clifford algebra. Spin structures are classified by the cohomology group $H^1(\check{X}_{\bullet}, \mathbb{Z}_2)$. If $h$ is a cocycle in $H^1(\check{X}_{\bullet}, \mathbb{Z}_2)$, then the spinor bundles associated with the different spin structures are reconstructed from the cocycle $\hat{g}h$ in $H^1(\check{X}_{\bullet}, \tx{Spin}_n)$ obtained by a pointwise product of locally defined maps getting values in $\text{Spin}_n$ (the elements of $\mathbb{Z}_2$ are embedded to $\text{Spin}_n$). 

The space of smooth sections of this bundle carries an algebraic structure which is given pointwise by the Clifford product. The product restricts to the space of invariant sections.  Suppose that a spin structure is fixed and let $F_{\Sigma}$ denote a spinor bundle.  The algebra of sections of the Clifford bundle acts on the sections of the spinor bundle. Since the groupoid action on the Clifford sections is by adjugation, the action of the Clifford sections on spinors can be restricted to define the action
\begin{eqnarray*}
\Gamma(\text{CL}( X_{\bullet}))^{\Theta} \times \Gamma(F_{\Sigma})^{\Theta} \rightarrow \Gamma(F_{\Sigma})^{\Theta}. 
\end{eqnarray*}
Consider invariant connections $\nabla_{\text{CL}}$ in $\text{CL}(X_{\bullet})$ and $\nabla$ in $F_{\Sigma}$. The  connection $\nabla$ is called Clifford compatible if it is compatible with the module structure: 
\begin{eqnarray*}
\nabla(\eta \cdot \psi) = (\nabla_{\text{CL}} \eta) \cdot \psi + \eta \cdot \nabla \psi,
\end{eqnarray*}
for all $\eta \in \Gamma(\text{CL}(X_{\bullet}))$ and $\psi \in \Gamma(F_{\Sigma})$. A complex Dirac bundle is a complex spinor bundle which is equipped with a unitary, Clifford compatible and $\Theta$-invariant groupoid connection and an inner product so that the unit vector fields $\gamma(u) \in \Gamma(\text{CL}(X_{\bullet}))$ act on the spinors as unitary transformations. We have the following result from \cite{Har14}. 
                      
\3 \noindent \textbf{Proposition 7.} A complex Dirac bundle exists on a spin proper {\' e}tale Lie groupoid.\3

Suppose that $X_{\bullet}$ is a spin proper {\' e}tale Lie groupoid. Choose $n$ linearly independent vector fields $e_i$ and let $e_i^*$ be the dual vector fields with respect to a fixed riemannian structure. The Dirac operator on a complex Dirac bundle is defined by 
\begin{eqnarray*}
\eth = \sum_{i = 1}^n \gamma(e_i^*) \nabla_{e_i}. 
\end{eqnarray*}
This operator is $\Theta$-invariant. In particular, it can be restricted to the operator 
\begin{eqnarray*}
\eth: \Gamma(F_{\Sigma})^{\Theta} \rightarrow \Gamma(F_{\Sigma})^{\Theta}. 
\end{eqnarray*}
Moreover, if $\psi,\psi' \in \Gamma(F_{\Sigma})^{\Theta}$, then 
\begin{eqnarray}\label{div}
(\eth \psi, \psi') - (\psi, \eth \psi') = \text{div}(V)  
\end{eqnarray}
where $V$ is a vector field such that $(V,W) = - (\psi, \gamma(W) \psi')$ for all $W \in \mathfrak{X}(X)$. The divergence $\text{div}(V)$ is a $\Theta$-invariant function. The goal is to define an essentially self-adjoint operator and therefore we need an inner product in the space of invariant sections of $F_{\Sigma}$ which kills the divergences of vector fields. We shall return to this Section 5. \3

\noindent \textbf{4.2.} Let $X_{\bullet}$ be a spin proper {\' e}tale groupoid. Fix a complex Dirac bundle $F_{\Sigma}$ and a Dirac operator $\eth$ on $X_{\bullet}$. Let $Y_{\bullet}$ be a representative of $X_{\bullet}$ and fix a Morita equivalence $\phi$. Then we have the induced bundles $\phi_{\#} \text{CL}(X_{\bullet})$ and $\phi_{\#} F_{\Sigma}$ on $Y_{\bullet}$. Let us write 
\begin{eqnarray*}
\phi_{\#} (\gamma(e))  =  \gamma_{\#} (e)  
\end{eqnarray*}
if $\gamma(e)$ is a section of $\text{CL}(X_{\bullet})$ to simplify notation.  Let  $\gamma_{\#}(e_1), \gamma_{\#}(e_2) \in \Gamma(\phi_{\#} \text{CL}(X_{\bullet}))^{\Xi}$. At $y \in Y$ the sections $\gamma_{\#}(e_1)$ and $\gamma_{\#}(e_2)$ get values in an equivalence class of fibres of $\text{CL}(X_{\bullet})$, as defined in 2.2. The bundle $\text{CL}(X_{\bullet})$ carries a pointwise Clifford product structure and when applied with the induced bundle, it gives the following pointwise multiplication in $\Gamma(\phi_{\#} \text{CL}(X_{\bullet}))^{\Xi}$: 
\begin{eqnarray*}
\gamma_{\#}(e_1) \cdot \gamma_{\#}(e_2) = \phi_{\#}(\gamma(e_1)\cdot \gamma(e_2))
\end{eqnarray*}
for all elements in $\Gamma(\phi_{\#} \text{CL}(X_{\bullet}))^{\Xi}$. It follows that the sections are subject to the anticommutation rules
\begin{eqnarray*}
\{ \gamma_{\#}(e_1), \gamma_{\#}(e_2) \} = (\gamma_{\#}(e_1), \gamma_{\#}(e_2))_{\#}. 
\end{eqnarray*}
Similarly, the pointwise Clifford module structure in $F_{\Sigma}$ induces a module structure in $\Gamma(\phi_{\#} F_{\Sigma})^{\Xi}$ under the action of  $\Gamma(\phi_{\#} \text{CL}(X_{\bullet}))^{\Xi}$ on $Y$. The action is given by 
\begin{eqnarray}\label{iclifford}
\gamma_{\#}(e) \cdot \phi_{\#} \psi = \phi_{\#}(\gamma(e) \cdot \psi). 
\end{eqnarray} 
In particular, as operators on the invariant spinor sections, the invariant Clifford sections satisfy 
\begin{eqnarray*}
\phi_{\#}(\gamma(e_1) \cdots \gamma(e_k)) = \gamma_{\#}(e_1) \cdots \gamma_{\#}(e_k)  = \phi_{\#} \circ \gamma(e_1) \cdots \gamma(e_k) \circ \phi_{\#}^{-1}
\end{eqnarray*}
for all $\gamma(e_1) \cdots \gamma(e_k) \in \Gamma(\text{CL}(X_{\bullet}))^{\Theta}$. Suppose that $u$ is a unit vector field on $X$. Then the following holds for the induced inner product
\begin{eqnarray*}
(\gamma_{\#}(u) \cdot \phi_{\#}\psi_1, \gamma_{\#} (u) \cdot \phi_{\#} \psi_2)_{\#} = ( \phi_{\#}\psi_1, \phi_{\#} \psi_2)_{\#}
\end{eqnarray*}
for all $\phi_{\#} \psi_1, \phi_{\#}\psi_2 \in \Gamma(\phi_{\#} \text{CL})^{\Xi}$. 

If $\nabla_{\text{CL}}$ is a Clifford connection in $\text{CL}(X_{\bullet})$ and $\nabla$ is a Clifford compatible connection in $F_{\Sigma}$, then the induced connection $\phi_{\#} \nabla$ satisfies the Clifford compatibility with the induced Clifford fields
\begin{eqnarray*}
\phi_{\#}\nabla (\phi_{\#}(\eta) \cdot \phi_{\#}(\psi)) =  \phi_{\#}\nabla_{\text{CL}}(\eta) \cdot \phi_{\#}(\psi) + \phi_{\#}(\eta) \cdot \phi_{\#} \nabla(\psi). 
\end{eqnarray*} 
for all $\phi_{\#}(\eta)\in \Gamma(\phi_{\#} \text{CL}(X_{\bullet}))^{\Xi}$ and $\phi_{\#} (\psi) \in \Gamma(\phi_{\#} F_{\Sigma})^{\Xi}$. 

Let us choose $n$ linearly independent vector fields, $e_i$, on the riemannian manifold $X$. Denote by $e_i^*$ the dual vector fields. The induced Dirac operator on $Y_{\bullet}$ is defined by: 
\begin{eqnarray*}
\phi_{\#} \eth &=& \sum_{i = 1}^n \gamma_{\#}(e_i^*) (\phi_{\#} \nabla)_{\phi_{\#} e_i}. 
\end{eqnarray*}
In general, $e_i$ and $e_i^*$ cannot all be invariant vector fields and so the induced fields $\phi_{\#} e_i$ and $\gamma_{\#}(e_i^*)$ are not well defined independently. However, the tensor field $\sum_i e_i^* \otimes e_i$ is an invariant section of a tensor bundle, and therefore the formula is well defined. \3

\noindent \textbf{Proposition 8.} The induced Dirac operator $\phi_{\#}\eth$ is an unbounded operator
\begin{eqnarray*}
\phi_{\#} \eth: \Gamma(\phi_{\#}F_{\Sigma})^{\Xi} \rightarrow \Gamma(\phi_{\#}F_{\Sigma})^{\Xi}
\end{eqnarray*}
and satisfies 
\begin{eqnarray}\label{divergence}
((\phi_{\#} \eth) \phi_{\#} \psi, \phi_{\#} \psi') - (\phi_{\#} \psi, (\phi_{\#} \eth) \phi_{\#} \psi' ) = \phi_{\#} \text{div}(V )  
\end{eqnarray}
If $V$ is a vector field on $X$ defined by $(V,W) = - (\psi, \gamma(W) \psi')$ for all $W \in \mathfrak{X}(X)$. \3

\noindent Proof. The statements in the proposition follow from 
\begin{eqnarray}\label{dirac}
\eth_{\#} = \phi_{\#} \circ \eth \circ \phi_{\#}^{-1}.
\end{eqnarray}
To verify this formula, let $\psi \in \Gamma(F_{\Sigma})^{\Theta}$, then 
\begin{eqnarray*}
\phi_{\#} (\gamma(e_i^*)\nabla_{e_i} \psi) &=&  \gamma_{\#}(e_i^*) \phi_{\#} (\nabla_{e_i}  \psi) \\
&=& \gamma_{\#}(e_i^*) (\phi_{\#} \nabla)_{\phi_{\#}(e_i)} (\phi_{\#} \psi) \\
&=&  \eth_{\#} (\phi_{\#} \psi). 
\end{eqnarray*}
In the first equality we use \eqref{iclifford}. The second equality follows from \eqref{connection}.  The induced operator is unbounded since $\eth$ is unbounded and $\phi_{\#}$ is an isomorphism. The equation \eqref{divergence} follows from the properties of the induced inner product, recall 3.2, and from the formulas \eqref{dirac} and \eqref{div}.\5 $\square$ \3

Suppose that the dimension of $X$ is even. Then there is an invariant chirality operator $\omega$ in $\Gamma(\text{CL}(X_{\bullet}))^{\Theta}$ and the induced chirality operator 
\begin{eqnarray*}
 \omega_{\#}: = \phi_{\#} \omega = \phi_{\#} \circ \omega \circ \phi_{\#}^{-1}. 
\end{eqnarray*}
The following  are direct consequences of $\omega^2 = 1$, \eqref{dirac} and Proposition 4:
\begin{eqnarray*}
\omega_{\#}^2 = 1, \5 \{\omega_{\#}, \eth_{\#} \} = 0,\5 [\omega_{\#}, \phi_{\#}(f)] = 0 
\end{eqnarray*}
for all invariant functions $\phi_{\#}(f) \in C^{\infty}(Y)^{\Xi}$. \3

\noindent \textbf{4.3.} Now let $X_{\bullet}$ and $Y_{\bullet}$ be both proper {\' e}tale groupoids and Morita equivalent to each other. Fix a Morita equivalence $\phi = (\varrho, Q, \alpha)$. We also assume that $X_{\bullet}$ and $Y_{\bullet}$ are both compact and spin, however, the property of being spin is formulated as an obstruction in a Lie groupoid cohomology and therefore it is sufficient to know that one of them is spin.  Now the theory of 4.1 can be applied with both groupoids $X_{\bullet}$ and $Y_{\bullet}$ but the next theorem proves that it is sufficient to fix any proper {\' e}tale groupoid in a Morita equivalence class to represent an orbifold. \3

\noindent \textbf{Theorem 1.}  Let $X_{\bullet}$ be a compact spin proper {\' e}tale groupoid and $\phi$ an {\' e}tale structure preserving Morita equivalence between $X_{\bullet}$ and $Y_{\bullet}$. Then $\phi$ induces a one to one correspondence between the complex Dirac bundles and Dirac operators on $X_{\bullet}$ and on $Y_{\bullet}$. \3

\noindent Proof. By definition, a Dirac bundle on $Y_{\bullet}$ (in the terminology of 4.1) is constructed by applying a spin lift to the structure cocycle of the tangent bundle $\tau Y_{\bullet}$. Then the lifted cocycle is composed with an element of $H^1(Y_{\bullet}, \mathbb{Z}_2)$ which fixes the spin structure. We first prove that the structure cocycle of $\tau Y_{\bullet}$ can be identified with the induced structure cocycle of $\phi_{\#} \tau X_{\bullet}$.  Let us fix the covers $\{U_a\}$ and $\{V_i\}$,  the {\v C}ech groupoids $\check{X}_{\bullet}$ and $\check{Y}_{\bullet}$ and the {\v C}ech bitorsor $\check{Q}$ as in Proposition 3. Since both $X_{\bullet}$ and $Y_{\bullet}$ are proper and {\' e}tale groupoids, the maps $\varrho$ and $\alpha$ are both local diffeomoprhisms. The locally defined sections $\beta_a^i$ can be assumed to be  diffeomorphisms. Associated to any arrow $\tau: y \rightarrow y'$ in $\Xi_{ij}$ there is the arrow 
\begin{eqnarray*}
\sigma: \varrho \circ \beta_b^j(y) \rightarrow \varrho \circ \beta_a^i(y') 
\end{eqnarray*}
in $\Theta_{ab}$ and the induced cocycle $\check{\phi}_{\#} g$ has the value $g_{ab}(\sigma) = (d \varphi_{\sigma})_{ab}$ at $\tau$. The differential $(d \varphi_{\sigma})_{ab}$ is written in the coordinates of $X_{\bullet}$ but we can use the locally defined diffeomorphisms $\varrho \circ \beta_a^i$ to write the components of the cocycle in the coordinates of $Y$. This gives a cocycle that is equivalent to $\phi_{\#} \check{g}$ which is given at $\tau \in \Xi_{ij}$ by 
\begin{eqnarray*}
(d\varrho \circ \beta_a^i)^{-1} (d\varphi_{\sigma})_{ab} (d\varrho \circ \beta_{b}^{j}) = d((\varrho \circ \beta_a^i)^{-1} \circ \varphi_{\sigma} \circ \varrho \circ \beta_b^j). 
\end{eqnarray*} 
The function $(\varrho \circ \beta_a^i)^{-1} \circ \varphi_{\sigma} \circ \varrho \circ \beta_b^j$ is a local diffeomoprhism arising from the action of $\Xi$ on $Y$ sending $y$ to $y'$. Now $Y_{\bullet}$ is {\' e}tale and the germs of such local diffeomoprhisms are unique, and so the value of the differential at $y \in V_j$ is equal to $(d\varphi_{\tau})_{ij}$ where $\varphi_{\tau}$ is a local diffeomorphism associated with the arrow $\tau$. It follows that the induced bundle $\phi_{\#} \tau X_{\bullet}$ is equivalent to the bundle $\tau Y_{\bullet}$. The Morita equivalence $\phi_{\#}$ is an isomorphism in cohomology and therefore sends the spin structures on $X_{\bullet}$ to the spin structures on $Y_{\bullet}$. Therefore the Clifford and spinor bundles on $Y_{\bullet}$ can be identified with induced Clifford and spinor bundles through the Morita equivalence $\phi$.  The map $\phi_{\#}$ induces a bijective correspondence between the inner product structures in vector bundles and spin connections, recall 3.1 and 3.2. It thus follows that the Dirac bundles on $Y_{\bullet}$ are precisely the induced Dirac bundles from $X_{\bullet}$ through any Morita equivalence. 

Let $e_i : 1 \leq i \leq n$ be a set of linearly independent vector fields on $X_{\bullet}$. The tensor field $\sum_i e_i^* \otimes e_i$ is $\Theta$-invariant and it can be written in the coordinates of $Y$ by applying the local diffeomorphisms $(d\varrho \circ \beta_a^i)^{-1}$ as above. This results in an invariant tensor field $\sum_i f_i^* \otimes f_i$ on $Y$ where $f_i$ are linearly independent at each point on $Y$ and the dual fields $f_i^*$ are defined with respect to the induced riemannian structure. So, under the identification of $\phi_{\#} F_{\Sigma}$ with a Dirac bundle on $Y_{\bullet}$, the induced Dirac operator can be written by $\sum_i \gamma(f_i^*) \nabla_{f_i}$ which is a Dirac operator on the Dirac bundle on $Y_{\bullet}$. \5 $\square$ \3

Let $X_{\bullet}$ and $Y_{\bullet}$ be the groupoids of Theorem 1. The structure cocycles in the multilinear bundles such as the dual $\tau X^*_{\bullet}$ and the exterior product bundles $\wedge^k \tau X_{\bullet}$ and $\wedge^k \tau X^*_{\bullet}$ are determined by the structure cocycle in $\tau X_{\bullet}$. It thus follows from the proof of Theorem 1 that the induced bundles $\phi_{\#} \tau X^*_{\bullet}$, $\phi_{\#} \wedge^k \tau X_{\bullet}$ and $\phi_{\#} \wedge^k \tau X^*_{\bullet}$ are isomorphic to the bundles $\tau Y^*_{\bullet}$, $\wedge^k \tau Y_{\bullet}$ and $\wedge^k \tau Y^*_{\bullet}$ on $Y_{\bullet}$. Now Proposition 4 implies the following. \3

\noindent \textbf{Corollary 2.} Let $X_{\bullet}$ be a proper {\' e}tale groupoids and $\phi$ an {\' e}tale structure preserving Morita equivalence between $X_{\bullet}$ and $Y_{\bullet}$. Then the linear map $\phi_{\#}$ of Proposition 5 is an isomorphism of vector spaces
\begin{eqnarray*}
\phi_{\#}: \Lambda^*(X_{\bullet})^{\Theta} \rightarrow \Lambda^*(Y_{\bullet})^{\Xi}.
\end{eqnarray*}

\section{Spectral Triples over Invariant Subalgebras}

The remaining task in the construction of a spectral geometric model over the invariant subalgebras of orbifolds is to develop an integration theory which can be applied for a Hilbert space completion of the space of invariant sections in a Dirac bundle. Since the orbit space of a groupoid is Morita invariant, it is natural to apply  an integration defined on the orbit space. Since the orbit space of a proper {\' e}tale groupoid has a structure of an orbifold, we can apply the well known integration theory on orbifolds. Then, in 5.2 the integration theory is extended over the Morita equivalence class. \3

\noindent \textbf{5.1.} Let $X_{\bullet}$ be a compact orientable proper {\' e}tale groupoid. The base of the groupoid $X$ can be given an open cover, $\{U_a\}$, such that the localized groupoids $\Theta^{U_a}_{U_a} \rightrightarrows U_a$ can be identified with action groupoids $G_a^a \ltimes U_a \rightrightarrows U_a$ for some finite subgroup $G_a^a$ of $GL_m$. The sets $U_a$ can be chosen so that the orbit space $|\Theta|$ gets a structure of an orbifold with respect to the charts $U_a$ and the projections $U_a \rightarrow |U_a|$ are determined by the restriction of the groupoid projection on $U_a$, \cite{MM03},\cite{MOE}. Then locally, on each $U_a$, one can define integration in the following way. $U_a$ has a stratification by orbit types. Let $U_a^{0}$ be the principal stratum and denote by $k_a$ the rank of the isotropy group at any point in the principal stratum (the rank is constant on strata). The integration of an invariant function $f$ on $X$ over $|U_a|$ is defined by 
\begin{eqnarray*}
\int_{|U_a|} f := \frac{k_a}{ |G_a^a|}\int_{U_a} f \nu
\end{eqnarray*}
where $\nu$ is the volume form associated with the groupoid riemannian structure. The coefficient in the integration is chosen because the projection map to the base restricts to a smooth covering map onto the principal stratum and there are $|G_a^a| / k $ sheets. On the other hand, the complement of the principal stratum has codimension greater or equal to one. The integration is extended to the orbit space $|\Theta|$ with a partition of unity which is a set of $G^a_a$-invariant smooth functions $\rho_a$ on each $U_a$ and these functions induce an ordinary topological partition of unity for $|\Theta|$ that is subordinate to $\{|U_a|\}$. Then the integral of an invariant function $f$ is defined by 
\begin{eqnarray*}
\int_{|\Theta|} f =  \sum_a \frac{k_a}{ |G_a^a|}\int_{U_a} \rho_a f \nu
\end{eqnarray*}
Now we can introduce the inner product in the space of $\Theta$-invariant spinor fields by 
\begin{eqnarray*}
\la \psi_1, \psi_2 \ra = \int_{|\Theta|} (\psi_1, \psi_2). 
\end{eqnarray*}
Let us denote by $L^2(F_{\Sigma})^{\Theta}$ the completion of the space $\Gamma(F_{\Sigma})^{\Theta}$ with respect to the norm associated with this inner product. The Stokes formula holds for the orbifold integration. In particular the Dirac operator can be extended to a self-adjoint unbounded operator on $L^2(F_{\Sigma})^{\Theta}$. 

\3 \noindent \textbf{5.2.} Let $X_{\bullet}$ be a compact spin proper {\' e}tale groupoid and $C^{\infty}(X)^{\Theta}$ the algebra of complex valued smooth invariant functions. Let $F_{\Sigma}$ be a spinor bundle. If $Y_{\bullet}$ is a representative of $X_{\bullet}$ through the Morita equivalence $\phi = (\varrho, Q, \alpha)$, then the induced space of invariant spinors $\phi_{\#} \Gamma( F_{\Sigma})^{\Xi}$ can be equipped with the obvious inner product so that the map $\phi_{\#}: \Gamma(F_{\Sigma})^{\Theta} \rightarrow \phi_{\#} \Gamma( F_{\Sigma})^{\Xi}$ extends to a unitary map of Hilbert spaces: for all $\psi_1, \psi_2 \in \Gamma(F_{\Sigma})^{\Theta}$ define 
\begin{eqnarray*}
\la \phi_{\#} \psi_1, \phi_{\#} \psi_2 \ra_{\#} &=& \int_{|\Theta|} \phi^{-1}_{\#} ( \phi_{\#}\psi_1, \phi_{\#}\psi_2 )_{\#} \\
&=& \int_{|\Theta|} (\psi_1, \psi_2)
\end{eqnarray*}
where $(\cdot, \cdot)$ is a $\Theta$-invariant pairing in $F_{\Sigma}$ and $(\cdot, \cdot)_{\#}$ is the induced pairing on $\phi_{\#} F_{\Sigma}$. Denote by $L^2(\phi_{\#} F_{\Sigma})^{\Xi}$ the Hilbert space completion.

\3 \noindent \textbf{Theorem 2.} Let $X_{\bullet}$ be a compact spin proper {\' e}tale groupoid. 
\begin{quote}
\textbf{1.} There is a finitely summable spectral triple $(C^{\infty}(X)^{\Theta}, \eth,L^2( F_{\Sigma})^{\Theta})$ which is equipped with the chirality element $\omega$ if the dimension of $X$ is even. 

\textbf{2.} If $\phi$ is a Morita equivalence between $X_{\bullet}$ and $Y_{\bullet}$, then there is a unitary equivalence of spectral triples 
\begin{eqnarray*}
(C^{\infty}(X)^{\Theta}, \eth,L^2( F_{\Sigma})^{\Theta}) \rightarrow (C^{\infty}(Y)^{\Xi}, \phi_{\#} \eth,L^2(\phi_{\#} F_{\Sigma})^{\Xi}) 
\end{eqnarray*}
and if the dimension of $X$ is even, this extends to a unitary equivalence of even spectral triples. 
\end{quote}

\3 \noindent Proof. The statement 1 is Theorem 1 of \cite{Har14}. 

The linear map $\Gamma(F_{\Sigma})^{\Theta} \rightarrow \Gamma(\phi_{\#} F_{\Sigma})^{\Xi}$ which sends $\psi$ to $\phi_{\#} \psi$ extends to a linear map of Hilbert spaces $U_{\phi}: L^2( F_{\Sigma})^{\Theta} \rightarrow  L^2(\phi_{\#} F_{\Sigma})^{\Xi} $ which is unitary with respect to the induced inner product by construction. It follows from Proposition 4 that as an operator on $\Gamma(\phi_{\#} F_{\Sigma})^{\Xi}$, the pointwise multiplication of the function $\phi_{\#}(f) \in C^{\infty}(Y)^{\Xi}$ extends to define the operator $U_{\phi} f U_{\phi}^{-1}$ on the Hilbert space $L^2(\phi_{\#} F_{\Sigma})^{\Xi}$.  As an operator on $L^2(\phi_{\#} F_{\Sigma})^{\Xi}$, the induced Dirac operator satisfies 
\begin{eqnarray*}
\phi_{\#} \eth = U_{\phi} \eth U_{\phi}^{-1} 
\end{eqnarray*}
and the same holds for the chirality operator. Therefore  $U_{\phi}$ defines a (even) unitary equivalence of (even) spectral triples. \5 $\square$

\3 \noindent \textbf{Note 2.} An orbifold groupoid is a Morita equivalence class of a proper {\' e}tale Lie groupoid. In this terminology, the content of the theorem is that one can associate a unitary equivalence class of invariant spectral triples with an orbifold. The choice of the proper {\' e}tale groupoid on which the spin geometric constructions are  based on is not unique. However, by Theorem 1, the same unitary equivalence class is obtained independently on this choice. 

\3 \noindent \textbf{Note 3.} On the level of spectral triples, the choice of the Morita equivalence $\phi$ has effect only on the description of the unitary map $U_{\phi}$, recall 1.6 and 2.2.

\section{Spectral Triples over Convolution Algebras}

\noindent \textbf{6.1.} Let $X_{\bullet}$ be a proper {\' e}tale groupoid and $C_c^{\infty}(\Theta)$ the algebra of compactly supported smooth complex valued functions equipped with the convolution product. If $\xi$ is a vector bundle on $X_{\bullet}$ then the convolution algebra $C_c^{\infty}(\Theta)$ can be represented on the space of compactly supported smooth sections in $\xi$ by setting
\begin{eqnarray*}
(f \cdot \psi)_x = \int_{\sigma \in \Theta^x} f(\sigma) (\varphi_{\sigma^{-1}}^{\#} \psi)_x   \mu^x(d \sigma) 
\end{eqnarray*}
for all $f \in C_c^{\infty}(\Theta)$ and $\psi \in \Gamma_c(\xi)$. We have used a left invariant Haar system $\mu$ so that $\mu^x$ is a left invariant measure in the fibres of $\Theta^x$ for all $x \in X$. 
 
\3 \noindent \textbf{Proposition 9.} Let $X_{\bullet}$ be a proper {\' e}tale groupoid. If $\xi$ is a vector bundle on $X_{\bullet}$, then the representation $\Gamma_c(\xi)$ of $C^{\infty}_c(\Theta)$ is faithful if and only if $X_{\bullet}$ is effective.  

\3 \noindent Proof. Suppose that the effectiveness does not hold. There is an open neighborhood $U$ in $X$, a pair of arrows $\sigma_1 \neq \sigma_2$ in $\Xi$ and local diffeomorphisms $\varphi_{\sigma^{-1}_1}$ and $\varphi_{\sigma^{-1}_2}$ associated with $\sigma_1^{-1}$ and $\sigma_2^{-1}$ which define the same map in $U$. Concretely, the source map has local bisections $\hat{\sigma_1}^{-1}: U \rightarrow W_1$ and $\hat{\sigma_2}^{-1}: U \rightarrow W_2$ so that $W_1$ and $W_2$ are local neighborhoods of $\sigma_1^{-1}$ and $\sigma_2^{-1}$ in $\Theta$, and then $\varphi_{\sigma_i^{-1}} : U \rightarrow t(W_i)$ is the composition $t \circ \hat{\sigma_i}^{-1}$ for $i = 1,2$. The local bisections are local diffeomorphisms because the groupoid is {\' e}tale. By choosing $U$ small enough we can assume that $W_1 \cap W_2 = \emptyset$.  If $f_1, f_2 \in C_c^{\infty}(\Theta)$ is a pair of functions supported in $W_1$ and $W_2$, then $\text{supp}(f_i) \cap \Theta^x$ is the point of the singlet set $W_i \cap \Theta^x$ for all $x \in U$ and $i = 1,2$. So, if $\psi \in \Gamma_c(\xi)$, then $f_1 \cdot \psi$ and $f_2 \cdot \psi$ have supports in $U$ and for all $x \in U$, 
\begin{eqnarray*}
(f_i \cdot \psi)_x &=& \int_{\kappa \in \Theta^x} f_i(\sigma) (\varphi^{\#}_{\kappa_i^{-1}} \psi)_y  \mu^y(d \kappa)  \\
&=& f_i(\hat{\sigma_i}^{-1}(x))(\varphi^{\#}_{\sigma_i^{-1}} \psi)_x  \lambda_i(x)
\end{eqnarray*}
for some strictly positive smooth functions $\lambda_i: U \rightarrow \mathbb{R}$ determined by the Haar measure. So, let us choose the values of $f_1$ and $f_2$ so that $f_1 \circ \hat{\sigma_1}^{-1} = f_2 \circ \hat{\sigma_2}^{-1}$. Then
\begin{eqnarray*}
(f_1 \cdot \psi)_x = (\Big(f_2 \frac{\lambda_1 \circ t}{\lambda_2 \circ t} \Big) \cdot \psi)_x
\end{eqnarray*}
for all $x \in X$ and for any $\psi \in \Gamma_c(\xi)$. The functions in the above formula are different since they have different supports in $\Theta$. So, the representation of $C^{\infty}_c(\Theta)$ is not faithful. 

Suppose that $X_{\bullet}$ is effective. If $U$ is any open set in $X$, then there is a point $x \in U$ so that the stabilizer at $x$ is the trivial group. On the other hand, if $z \in U$ is any point, then there is a neighborhood $N_{z}$ of $z$ so that the groupoid $X_{\bullet}$ localizes to an action groupoid $\Theta_z^z \ltimes N_z \rightrightarrows N_z$ over $N_z$. So, if we take $x = z$, then we find that there is a neighborhood $N_x$ of $x$ so that the set of arrows $s^{-1}(N_x) \cap t^{-1}(N_x)$ consist of the identity arrows of $N_x$ only. Suppose that $\sigma: x' \rightarrow x$. We can take $N_x$ to be arbitrarily small, and for this reason we can assign the diffeomorphism $\varphi_{\sigma^{-1}}: N_x \rightarrow \varphi_{\sigma^{-1}}(N_x)$ with $\sigma^{-1}$. Observe that the arrows $s^{-1}(N_x) \cap t^{-1}(\varphi_{\sigma}^{-1} N_x)$ are exactly the arrows which are determined by $\varphi_{\sigma^{-1}}$ in $N_x$. Namely, if $\kappa: x_1 \rightarrow x_2$ is an arrow $N_x \rightarrow  \varphi_{\sigma^{-1}}(N_x)$ then $\varphi_{\sigma} \circ \kappa$ would define an arrow in $N_x$ which needs to be the unit arrow. Thus, $\varphi_{\sigma}^{-1}(x_1) = x_2$. 

Let $f_1, f_2 \in C_c^{\infty}(\Theta)$ so that $f_1 \neq f_2$. Then there is a neighborhood $W \subset \Theta$ so that $f(\sigma) \neq g(\sigma)$ at every point $\sigma \in W$ and we can choose $W$ to be small enough so that $t$ restricts to a diffeomoprhism in $W$. Now $t(W)$ is an open subset in $X$ and so the discussion of the previous paragraph applies. Let us fix $x \in t(W)$  so that $x$ has trivial isotropy. There is the arrow $\sigma \in W$ so that $t(\sigma) = x$. According to the analysis above, we can choose a neighborhood $N_x$ for $x$ so that there is the local diffeomorphism $\varphi_{\sigma^{-1}}: N_x \rightarrow \varphi_{\sigma^{-1}}(N_x)$ and $s^{-1}(N_x) \cap t^{-1}(\varphi_{\sigma^{-1}}(N_x)) \subset W$. Choose $\psi$ to be a nonzero element in $\Gamma_c(\xi)$ whose support lies in $\varphi_{\sigma^{-1}}(N_x)$. If $z \in N_x$ there will be contributions to the integral over the $\Theta^z$-fibre only from the element of the singlet $W  \cap \Theta^z$. For the values $z \in N_x$ where $\psi_{\varphi_{\sigma}^{-1}(z)}$ is nonzero, we have
\begin{eqnarray*}
f_1(\kappa) (\varphi^{\#}_{\sigma^{-1}} \psi)_z \neq f_2(\kappa) (\varphi^{\#}_{\sigma^{-1}} \psi)_z
\end{eqnarray*}
where $\kappa$ is the arrow in $ W \cap \Theta^z$. Thus, $f_1 \cdot \psi \neq f_2 \cdot \psi$. \5 $\square$

\3 \noindent \textbf{6.2.} Suppose that $X_{\bullet}$ is a spin proper {\' e}tale groupoid and $\phi$ is an {\' e}tale structure preserving Morita equivalence between $X_{\bullet}$ and $Y_{\bullet}$. 

The Dirac operator introduced in 4.1 is a $\Theta$-invariant differential operator, and for this reason the induced Dirac operator of 4.2 can be defined on the space of  smooth compactly supported sections $\Gamma_c(\phi_{\#} F_{\Sigma})$. Let $y \in Y$ and choose an open neighborhood $V \subset Y$ for $y$ which is small enough so that there is a local section $\beta: Y \rightarrow Q$ of $\alpha$, and the composition $\varrho \circ \beta$ restricts to a diffeomoprhism. Now we can fix a local representative for the fibres of $\phi_{\#}F_{\Sigma}$ over $V$: if $y \in V$, then $(\phi_{\#}F_{\Sigma})_y$ is represented by $(F_{\Sigma})_{\varrho \circ \beta(y)}$ and so $\phi_{\#}F_{\Sigma}$ on $V$ can be identified with the pullback bundle $(\varrho \circ \beta)^* F_{\Sigma}$. The sections of $\phi_{\#}F_{\Sigma}$ on $V$ are pullbacks of the sections on $F_{\Sigma}$ under the diffeomorphism $\varrho \circ \beta$. The Dirac operator at $y \in Y$ is represented by 
\begin{eqnarray*}
(\phi_{\#} \eth)_q(y) = (\varrho \circ \beta)^*  \circ \eth_{\varrho(q)} \circ ((\varrho \circ \beta)^{-1})^*
\end{eqnarray*}
if $\beta(y) = q$. The $\Xi$-invariance follows since the construction factors through $\varrho$. When restricted to the $\Xi$-invariant compactly supported sections, this reduces to the operator \eqref{dirac}. To see that this gives a well defined operator on $\Gamma_c(\phi_{\#} F_{\Sigma})$, we need to verify that the local descriptions associated with the different choices of the $\alpha$-fibre parameters are properly related to each other. This is the content of the next lemma. 

\3 \noindent \textbf{Lemma 3.} If $q$ and $\sigma \cdot q$ are in $\alpha^{-1}(y)$, then 
\begin{eqnarray*}
\rho^{-1}(\sigma) (\phi_{\#} \eth)_{\sigma \cdot q}(y) \rho(\sigma) = (\phi_{\#} \eth)_q(y).
\end{eqnarray*}

\noindent Proof. We have two different representatives of $\phi_{\#} \eth$ associated with two local sections $\beta_1$ and $\beta_2$ of $\alpha$. Suppose that $\beta_1(y) = q$ and $\beta_2(y) = \sigma \cdot q$ and that $V$ is their common domain. The values of $\beta_1$ and $\beta_2$ are connected by the action of $\Theta$ on $Q$ and so, by shrinking $V$ is necessary, we can write $\beta_2 = \hat{\varphi_{\sigma}} \circ \beta_1$ in the notation of Lemma 2. A straightforward application of Lemma 2.1 with the $\Theta$-invariance of $\eth$ gives the claim. \5 $\square$

\3 \noindent \textbf{Proposition 10.} Suppose that $\phi$ is an {\' e}tale structure preserving Morita equivalence between $X_{\bullet}$ and $Y_{\bullet}$. Then for all $g \in C^{\infty}(Y)$ and $y \in Y$ 
\begin{eqnarray*}
[\phi_{\#} \eth, g]_y = \sum_{i} (f_i g)_y \gamma(f_i^*)_y
\end{eqnarray*}
where $f_i$ are linearly independent vector fields on $Y$.  

\3 \noindent Proof. Locally, around $y \in Y$, the Dirac operator is represented by a pullback Dirac operator through a diffeomorphism. In the coordinates of $Y$, the pullback Dirac operator is an ordinary Dirac operator associated with the pullback metric and the pullback spin structure. So, if $g \in C^{\infty}(Y)$, the claim follows from the standard analysis, \cite{LM89}. \3  

In the induced spinor bundle  $\phi_{\#} F_{\Sigma}$ we have the inner product defined in 3.2 which gives the pointwise pairing of sections $\Gamma_c(\phi_{\#}F_{\Sigma}) \otimes \Gamma_c(\phi_{\#}F_{\Sigma}) \rightarrow C^{\infty}_c(Y)$ defined by
\begin{eqnarray*}
(\psi_1,\psi_2)_{\#,y} = (\psi_1, \psi_2)_{\varrho(q)}
\end{eqnarray*}
where $q$ is an arbitrary point in the $\alpha$-fibre over $y$.  So, over $Y$ we have the $L^2$-inner product in the space of sections of the bundle $\phi_{\#} F_{\Sigma}$ which is defined by 
\begin{eqnarray*}
\la \psi_1, \psi_2 \ra_{\#} = \int_Y ( \psi_1, \psi_2)_{\#} 
\end{eqnarray*}
for all $\psi_1, \psi_2 \in \Gamma_c(\phi_{\#} F_{\Sigma})$. The Hilbert space completion is denoted by $L^2(\phi_{\#} F_{\Sigma})$. 

\3 \noindent \textbf{Theorem 3.} Let $X_{\bullet}$ be an effective spin proper {\' e}tale groupoid over a complete base manifold $X$.
\begin{quote}
\textbf{1.} There is a finitely summable spectral triple $(C^{\infty}_c(\Theta),\eth, L^2(F_{\Sigma}))$ which is equipped with the chirality operator $\omega$ if the dimension of $X$ is even. 

\textbf{2.} If $\phi$ is an {\' e}tale structure preserving Morita equivalence between $X_{\bullet}$ and $Y_{\bullet}$ and if the base manifold $Y$ is complete, then there is a finite summable spectral triple 
\begin{eqnarray*}
(C^{\infty}_c(\Xi), \phi_{\#} \eth, L^2(\phi_{\#}F_{\Sigma}))
\end{eqnarray*}
which is equipped with the chirality operator $\omega_{\#}$ if the dimension of $X$ is even. 
\end{quote}

\noindent Proof. The Dirac operator $\eth$ is the $\Theta$-invariant Dirac operator defined in 4.1. The statement 1 is Theorem 2 of \cite{Har14}. 

The effectiveness of a groupoid is preserved under an {\' e}tale structure preserving Morita equivalence and so $X_{\bullet}$ is effective if and only if $Y_{\bullet}$ is effective, \cite{MM03}. It follows from Proposition 9 that the induced representation is faithful. 

To see that the Dirac operator has bounded commutators with the representation operators of $C^{\infty}_c(\Xi)$ on $\Gamma_c(\phi_{\#} F_{\Sigma})$ we can approach as in Proposition 10. Namely, if $y \in Y$ we can choose a local neighborhood $V$ of $y$ so that the spinor bundle and the Dirac operator can be represented by a pullback Dirac bundle and a pullback Dirac operator through a diffeomorphism over $V$. So, locally the statement follows from the corresponding result in \cite{Har14}. The claim holds globally over $Y$ because $Y$ can be covered with the subsets which have the properties of $V$. Therefore the commutators $[\phi_{\#}\eth, f]$ extend to bounded operators on $L^2(\phi_{\#} F_{\Sigma})$ for all $f \in C^{\infty}_c(\Xi)$. The finite summability holds since $\phi_{\#} \eth$ is a first order elliptic differential operator on  $\Gamma_c(F_{\Sigma})$. Note that it is sufficient to check these conditions locally and they hold locally by the discussion above.

We need to show that the closure of $\phi_{\#} \eth$ in $L^2(\phi_{\#}F_{\Sigma})$ is a self-adjoint operator. The operator $\phi_{\#} \eth$ and its adjoint $(\phi_{\#} \eth)^*$ both have Hilbert space closures. As above, we use the property that in a small enough subsets of $V$, the operator $\phi_{\#} \eth$ can be represented by a locally defined Dirac operator on $Y$. So by the usual partition of unity argument, the formal self-adjointness of $\phi_{\#} \eth$ holds. The inclusion $\text{dom}(\phi_{\#}\eth) \subset \text{dom}((\phi_{\#} \eth)^*)$ follows. For the reverse inclusion, $\text{dom}((\phi_{\#}\eth)^*) \subset \text{dom}(\phi_{\#} \eth)$, the standard proof can be applied without any difficulties: \cite{LM89} Theorem 5.7. Since the base manifold is complete, and $\phi_{\#} \eth$ satisfies the formula in Proposition 10, we can reduce this problem to the case of spinors with compact support in $\text{dom}(\phi_{\#} \eth^*)$. Now we can use the partition of unity argument to reduce the problem to an arbitrarily small open domain in $Y$ where the statement holds by the usual Dirac operator analysis. 

The chirality section $\omega$ induces a $\Xi$-invariant section $\omega_{\#}$, recall 4.2,  which acts on the sections of $\phi_{\#} F_{\Sigma}$. Using the $\Xi$-invariance we find that for all $f \in C^{\infty}_c(\Xi)$ 
\begin{eqnarray*}
(\omega_{\#})_y (f \cdot \psi)_y &=& (\omega_{\#})_y \int_{\tau \in \Xi^y} f(\tau) (\varphi^{\#}_{\tau^{-1}} \psi)_y \mu^y (d \tau) \\
&=& \int_{\tau \in \Xi^y} f(\tau) \varphi^{\#}_{\tau^{-1}} ((\omega_{\#}) \psi)_y \mu^y (d \tau) \\
&=& f \cdot ((\omega_{\#}) \psi)_y 
\end{eqnarray*}
and so the action of $\omega_{\#}$ on $\Gamma_c(\phi_{\#} F_{\Sigma})$ commutes with the representation of the function algebra. Then $ \omega_{\#}$ extends to a bounded self-adjoint operator on $L^2(\phi_{\#}F_{\Sigma})$. It follows that $\omega_{\#}$ defines an even structure in the spectral triple. \5 $\square$

\section*{Appendix: On Morita Equivalence}

\noindent \textbf{A.1.} There is another and better known formulation of the Morita equivalence of Lie groupoids. I have chosen to follow the conventions in 1.1 because those are convenient in the study of geometric structures. Let $X_{\bullet}$ and $Y_{\bullet}$ be Lie groupoids. A morphism of Lie groupoids, $\theta_{\bullet}: X_{\bullet} \rightarrow Y_{\bullet}$, is called a weak equivalence if 
\begin{quote} 
	\textbf{1.} $\theta_{0}: X_{(0)} \rightarrow Y_{(0)}$ is a surjective submersion.
	
	\textbf{2.} The following diagram is cartesian: 
	\begin{center}
	\isonelio{X_{(1)}}{X_{(0)}\times X_{(0)}}{Y_{(1)}}{Y_{(0)} \times Y_{(0)}}{(s,t)}{(s,t)}{\theta_1}{\theta_0 \times \theta_0}  
	\end{center}
	\end{quote}

\noindent The following result is well known, \cite{BX11}. \3
	
\noindent \textbf{Proposition 11.} Let $X_{\bullet}$ and $Y_{\bullet}$ be Lie groupoids. The following are equivalent:
\begin{quote}
\textbf{1.} There is a Lie groupoid $Z_{\bullet}$ and a pair of weak equivalences: $X_{\bullet} \leftarrow Z_{\bullet} \rightarrow Y_{\bullet}$. 

\textbf{2.} The groupoids $X_{\bullet}$ and $Y_{\bullet}$ are Morita equivalent. 
\end{quote}

\noindent \textbf{A.2.} As is evident from the discussion of the Morita equivalence in 1.1, the orbit space of the groupoid is invariant under Morita equivalences. The following result clarifies how the base manifolds $X$, $Y$ and the isotropy are related under a Morita equivalence between $X_{\bullet}$ and $Y_{\bullet}$. \3

\noindent \textbf{Proposition 12.} Let $\phi = (\varrho, Q, \alpha)$ be a Morita equivalence between $X_{\bullet}$ and $Y_{\bullet}$. Let $y \in Y$ and let the rank of the isotropy group $\Xi_y^y$ at $y \in Y$ be equal to $k$. The following hold:
\begin{quote}
\textbf{1.} The fibre $\alpha^{-1}(y)$ has a partition into subsets that are invariant under the action of $\Xi_y^y$. The cardinality of each block is $k$. 

\textbf{2.} The map $\varrho$ restricts to a constant map on the blocks and if $\varrho(q) = x$ for some $q \in \alpha^{-1}(y)$, then $\Theta_x^x$ has rank $k$.
\end{quote}

\noindent Proof. If $y \in Y$ the fibre $\alpha^{-1}(y)$ is stable under the action of $\Xi_y^y$ and so the fibre has a partition into the orbits of the action. Each orbit must have cardinality $k$ since the action is free by definition. 

Fix some $q_1$ so that $\alpha(q_1) = y$. Write $\{q_i: 1 \leq i \leq k\}$ for the orbit under the $\Xi_y^y$ action. The elements of the isotropy group $\Xi_y^y$ act as arrows of $\Xi$ and so the action is in the direction of the $\varrho$-fibres. Write $\varrho(q_i) = x \in X$ which is independent on the choice $1 \leq i \leq k$. The points $\{q_i\}$ are also in the same $\alpha$-fibre because the action of $\Xi_y^y$ preserves the fibre at $y$. It follows that $\Theta_x^x$ acts freely and transitively on the set $\{q_i\}$.  Therefore the order of $\Theta_x^x$ is $k$. \5 $\square$ \3

Proposition 12 implies that for each order $k$ singularity in $Y$ there is an order $k$ singularity in $X$. However, it is not true in general that the singular points in $X$ and the singular points in $Y$ are in one to one correspondence. For example, there may be an isolated singularity in $X$ of order $k$ and in $Y$ the order $k$ singular points are localized along an orbit of  the $\Xi$-action. In a simplified case one can take $X_{\bullet}$ to be the groupoid $\mathbb{Z}_2 \rightrightarrows *$ (this is the groupoid theoretic model for the group $\mathbb{Z}_2$, $*$ is a point) and $Y_{\bullet}$ to be the action groupoid $\mathbb{T} \ltimes \mathbb{T} \rightrightarrows \mathbb{T}$ where the $\mathbb{T} = \mathbb{R} / 2 \pi \mathbb{Z}$ action on itself is given by $([x],[y]) \mapsto ([2x + y])$. We use the notation $[x] = x \text{ mod } 2 \pi$. The Morita bitorsor can be chosen by $Q = \mathbb{R} / 4 \pi \mathbb{Z}$. Then $\phi = (\varrho, Q, \alpha)$ is a Morita equivalence so that $\varrho: Q \rightarrow *$ is the constant map and $\alpha: Q \rightarrow \mathbb{T}$ is the map which $x \text{ mod } 4 \pi$ to $x \text{ mod } 2\pi$. The nontrivial element of $\mathbb{Z}_2$ acts on $Q$ by $2 \pi$-translations and the circle $\mathbb{T}$ acts on $Q$ by $([x], [q]_{4 \pi}) \mapsto [2x + q]_{4 \pi}$.


\begin{thebibliography}{000}
\bibitem{Beh04} Behrend K.: Cohomology of Stacks, in: Intersection theory and moduli, ICTP Lect. Notes, XIX, Abdus Salam Int. Cent. Theoret. Phys., Trieste, 249-294 (2004)
\bibitem{BX11} Behrend K., Xu P.: Differentiable Stacks and Gerbes, J. Symplectic Geom. 9, 285-341 (2011)
\bibitem{A1} Bratteli O., Elliott G.A., Evans D.E., Kishimoto A.: Noncommutative spheres. I, Internat. J. Math. 2, 139-166 (1991)
\bibitem{BF12}  Brzezinski T., Fairfax S. A.: Quantum Teardrops, Commun. Math. Phys. 316, 151-170 (2012)
\bibitem{A3} Brzezinski T., Zielinski B., Quantum principal bundles over quantum real projective spaces, J. Geom. Phys. 62, 1097-1107 (2012)
\bibitem{Con94} Connes A.: Noncommutative Geometry (Academic Press, San Diego, CA 1994) 
\bibitem{Con13}  Connes A.: On the spectral characterization of manifolds, J. Noncom. Geom. 7, 1-82 (2013)
\bibitem{GL13} Gorokhovsky A., Lott J.: The index of a transverse Dirac-type operator: the case of abelian Molino sheaf, J. Reine Angew. Math. 678, 125-162 (2013)
\bibitem{GH06} Guruprasad K., Haefliger A.: Closed Geodesics on Orbifolds, Topology 46, 611–641 (2006)
\bibitem{Har14} Harju A. J.: Spectral Triples on Proper {\' E}tale Groupoids (to appear in Journal of Noncommutative Geometry) arXiv:1402.6255 (2014)
\bibitem{Har14b} Harju A. J.: Dirac Operators on Quantum Weighted Projective Spaces, arXiv:1402.6251 (2014)
\bibitem{A2} Hong J.H., Szymanski W.: Quantum lens spaces and graph algebras, Pacific J. Math. 211, 249-263 (2003)
\bibitem{LTX07} Laurent-Gengoux C., Tu J-L., Xu P.: Chern-Weil map for principal bundles over groupoids, Math Z. 255, 451-491 (2007)
\bibitem{LM89} Lawson H. B., Michelsohn L-M.: Spin Geometry (Princeton University Press 1989)
\bibitem{Moe02} Moerdijk I.: Orbifolds as groupoids: an introduction. Orbifolds in mathematics and physics (Madison, WI, 2001), Contemp. Math. 310, 205-222 (Amer. Math. Soc., Providence, RI, 2002)
\bibitem{MM03} Moerdijk I.,  Mrcun J.: Introduction to Foliations and Lie Groupoids (Cambridge University Press 2003)
\bibitem{MOE} Moerdijk I., Pronk D. A.: Orbifolds, Sheaves and Groupoids, $K$-Theory 12, 3-21 (1997)
\bibitem{RV08} Rennie, A. C., Varilly, J. C.: Orbifolds are not commutative geometries, J. Australian Math. Soc. 84, 109-116 (2008)
\bibitem{Sat56} Satake I.: On a generalization of the notion of manifold,  Proc. Nat. Acad. Sci. U.S.A. 42, 359-363 (1956)
\bibitem{SV13} Sitarz A., Venselaar J. J.: Real spectral triples on 3-dimensional noncommutative lens spaces, arXiv:1312.5690 (2013)
\end{thebibliography}
\end{document}